\documentclass[letterpaper,12pt,peerreviewca,draftcls]{IEEEtran}

\usepackage{csmMBFeb07}

\title{%
	Sensitivity Analysis of Oscillator Models \\
	in the Space of Phase-Response Curves\\
	{\Large Oscillators as open systems}
}
\author{Pierre Sacr\'e and Rodolphe Sepulchre --- \today}

\usepackage{graphicx}

\newcommand{\myfigurename}{Figure}

\newcommand{\mytablename}{Table}

\usepackage{amsmath}		
\usepackage{amssymb}
\usepackage{amsthm}
\usepackage{mathtools}

\theoremstyle{plain}

\theoremstyle{definition} 
\newtheorem{defn}{Definition}

\theoremstyle{remark} 
\newtheorem*{rem}{Remark}

\newcommand{\field}[1]{\mathbb{#1}} 
 
\newcommand{\Rbb}{\field{R}}
\newcommand{\Nbb}{\field{N}}
\newcommand{\Sbb}{\field{S}}
\newcommand{\Zbb}{\field{Z}}

\newcommand{\Hilbert}{\mathcal{H}}
\newcommand{\Lp}[3]{\mathcal{L}_{#1}(#2,#3)}

\DeclareMathOperator{\dist}{dist}
\DeclareMathOperator{\Shift}{Shift}

\newcommand{\cconj}[1]{{#1}^{*}}
\newcommand{\trans}[1]{#1^T}

\DeclareMathOperator*{\argmin}{arg\,min}
\DeclareMathOperator*{\argmax}{arg\,max}
\DeclareMathOperator{\diag   }{diag}
\DeclareMathOperator{\grad   }{grad}
\DeclareMathOperator{\range  }{range}

\newcommand{\tgspace}[2]{T_{#1}{#2}}
\newcommand{\dirder}[3]{D#1(#2)[#3]}
\newcommand{\Euclgrad}[3]{\nabla_{#3}#1(#2)}

\DeclarePairedDelimiter{\abs}{\lvert}{\rvert}
\DeclarePairedDelimiterX{\innerp}[2]{\langle}{\rangle}{#1,#2}
\DeclarePairedDelimiter{\norm}{\lVert}{\rVert}

\newcommand{\eqdef}{\mathrel{:=}}
\newcommand{\reveqdef}{\mathrel{=:}}

\makeatletter
\newcommand{\prap}[1]{\allowbreak\if@display\mkern18mu\else\mkern8mu\fi(#1)}
\newcommand{\pwrap}[1]{\prap{{\operator@font wrap\;to}\mkern6mu#1}}
\makeatother

\usepackage{nomencl}  
\makenomenclature

\usepackage{color}
\usepackage{cite}
\usepackage{url}

\usepackage{enumerate}

\newcommand{\zeroinput}{\mathbf{0}}

\newcommand{\fvec}{f}
\newcommand{\hvec}{h}

\newcommand{\flow  }{\phi}
\newcommand{\perorb}{\gamma}
\newcommand{\persol}{x^{\perorb}}

\newcommand{\PRC  }{Q}
\newcommand{\PTC  }{R}
\newcommand{\iPRCu}{q}
\newcommand{\iPRCx}{p}

\newcommand{\PLC  }{\hat{\varphi}}

\newcommand{\mesh}{\Pi}

\newcommand{\nlambda}{\lambda_0}

\begin{document}

\maketitle
\CSMsetup


\clearpage

Oscillator models---whose steady-state behavior is periodic rather than constant---are fundamental to rhythmic modeling and they appear in many areas of engineering, physics, chemistry, and biology~\cite{Winfree:1980ue,Kuramoto:1984wo,Glass:1988ub,Goldbeter:1996uo,Pikovsky:2001et,Izhikevich:2007vr}.
Many oscillators are, by nature, open dynamical systems, that is, they interact with their environment~\cite{Sepulchre:2006vk}. Whether they function as clocks, information transmitters, or rhythm generators, these oscillators have the robust ability to respond to a particular input (entrainment) and to behave collectively in a network (synchronization or clustering).

The phase response curve of an oscillator has emerged as a fundamental input--output characteristic of oscillators~\cite{Winfree:1980ue}. Analogously to the static (zero-frequency) gain of a transfer function, the phase response curve measures a steady-sate (asymptotic) property of the system output in response to an impulse input. For the zero-frequency gain, the measured quantity is the integral of the response; for the phase response curve, the measured quantity is the  phase shift between the perturbed and  unperturbed responses. Because of the periodic nature of the steady-state behavior, the magnitude and the sign (advance or delay) of this phase shift depend on the phase of the impulse input. The phase response curve is therefore a curve rather than a scalar.
In many situations, the phase response curve can be determined experimentally and provides unique data for the systems analysis of the oscillator. Alternatively, numerical methods exist to compute the phase response curve from a mathematical model of the oscillator. 
The phase response curve is the fundamental mathematical information required to reduce an $n$-dimensional state-space model to a one-dimensional (phase) center manifold of a hyperbolic periodic orbit.

Motivated by the prevalence of the input--output representation in experiments and the growing interest in system-theoretic questions related to oscillators, this article extends fundamental concepts of systems theory to the space of phase response curves. 
Comparing systems with a proper metric has been central to systems theory (see \cite{Zames:1980up,ElSakkary:1985gr,Vinnicombe:2000wt,Georgiou:2007cw} for exemplative milestones). In a similar spirit, this article aims to endow the space of phase response curves with the right metrics (accounting for natural equivalence properties) and sensitivity analysis tools. This framework provides mathematical and numerical grounds for robustness analysis and system identification of oscillator models. 
Although classical in their definitions, several of these tools appear to be novel, particularly in the context of biological applications.

The focus of the article is on oscillator models in systems biology and neurodynamics---two areas where sensitivity analysis is particularly useful to assist the increasing focus on quantitative models.
In systems biology, phase response curves have been primarily studied in the context of circadian rhythms models~\cite{DeCoursey:1960cw,Goldbeter:1996uo,Goldbeter:2002iy}.
A circadian oscillator is at the core of most living organisms that need to adapt their physiological activity to the 24 hours environmental cycle associated with earth's rotation (for example variations in light or temperature condition).
This oscillatory system is capable of exhibiting oscillations with a period close to 24 hours in constant environmental condition and of locking its oscillations (in frequency and phase) to an environmental cue with a period equal to 24 hours.
In neurodynamics, the use of phase response curves is more recent but increasingly popular~\cite{Schultheiss:2012bz}.
A spiking oscillator is the repeated discharge of action potentials by a neuron, which is the basis for neural coding and information transfer in the brain. 
This oscillatory system is capable of exhibiting oscillations on a wide range of period---from $0.001$ to $10$ secondes---and of behaving collectively in a neural network.
Phase response curves are also used in many other areas of sciences and engineering (planar particle kinematics, Josephson junctions, alternating current power networks, etc.) for which the reader is referred to the abundant literature (see for example the pioneering contributions~\cite{Winfree:1967vf,Kuramoto:1975ki,Winfree:1980ue,Kuramoto:1984wo,Strogatz:2000wx,Strogatz:2003tm} and the detailed review~\cite[and references therein]{Dorfler:2013fk}).

The results of the article primarily draw out from the Ph.D. dissertation of the first author~\cite{Sacre:2013ys}.
A preliminary version of this work was presented in \cite{Sacre:2011vg}. The first case study on circadian rhythms was discussed in detail in \cite{Sacre:2012fk}.

The article is organized as follows.
``Phase Response Curves from Experimental Data'' presents the concept of phase response curves derived from phase-resetting experiments.
``Phase Response Curves from State-Space Models'' reviews the notion of phase response curves characterizing the input--output behavior of an oscillator model in the neighborhood of an exponentially stable periodic orbit.
``Metrics in the Space of Phase Response Curves'' defines several relevant metrics on (nonlinear) spaces of phase response curves induced by natural equivalence properties.
``Sensitivity Analysis in the Space of Phase Response Curves'' develops the sensitivity analysis for oscillators in terms of the sensitivity of its periodic orbit and its phase response curve.
``Applications to Biological Systems'' illustrates how these tools solve system-theoretic problems arising in biological systems, including robustness analysis, system identification, and model classification.

The main developments of the article are supplemented by several supporting discussions.
``A Brief History of Phase Response Curves'' sets the use of phase response curves in its historical context.
``Phase Maps'' defines the key ingredients for studying oscillator models on the unit circle.
``From Infinitesimal to Finite Phase Response Curves'' provides details on the mathematical relationship between finite and infinitesimal phase response curves.
``Basic Concepts of Differential Geometry on Manifolds'' and ``Basics Concepts of Local Sensitivity Analysis'' present distinct features of differential geometry and sensitivity analysis used in this article, respectively.
``Numerical Tools'' provides the numerical tools to turn the abstract developments into concrete algorithms.
The notation is defined in ``List of Symbols''.

\clearpage
\section{Phase Response Curves from Experimental Data} \label{sec:prc_exp}

The interest of a biologist in an oscillator model comes through the observation of a rhythm, that is, the regular repetition of a particular event. Examples include the onset of daily locomotor activity of rodents, the initiation of an action potential in neural or cardiac cells, or the onset of mitosis in cells growing in tissue culture 
(see ``A Brief History of Phase Response Curves'').
One of the simplest modeling experiments is to perturb the oscillatory behavior for a short (with respect to the oscillation period) duration and record the altered timing of subsequent repeats of the observable event. Once the system has recovered its prior rhythmicity, the phase of the oscillator is said to have reset.
In general, the phase reset depends not only on the perturbation itself (magnitude and shape) but also on its timing (or phase) during the cycle.
This section formalizes the basic experimental paradigm of phase-resetting experiments and describes the concept of phase response curves following the terminology in \cite{Winfree:1980ue} and \cite{Glass:1988ub}.

\begin{figure}[p]
	\centering
	\includegraphics{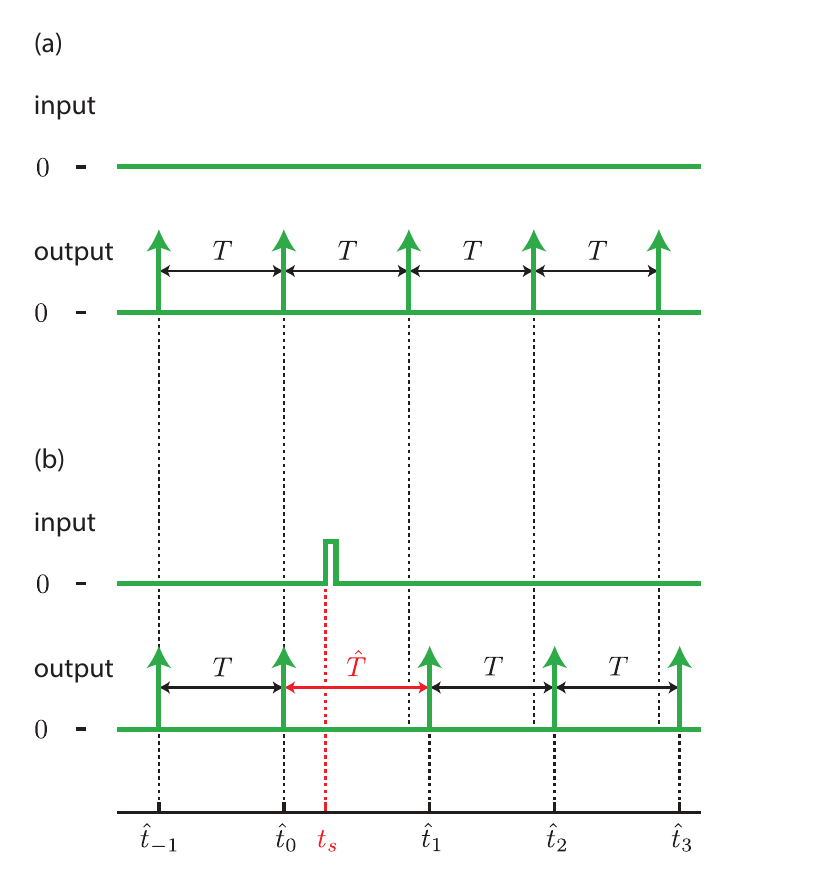}
	\caption{%
		Schematic representation of a phase-resetting experiment.
		(a)~In isolated conditions (closed system), the observable event (vertical arrow) occurs every  $T$ units of time. 
		(b)~Following a phase-resetting stimulus at time $(t_s - \hat{t}_0)$ after an event (open system), the successive observable event times~$\hat{t}_i$, for $i\in\Nbb_{>0}$, are altered.
		$\hat{T} \eqdef \hat{t}_1 - \hat{t}_0$ denotes the time interval between the pre- and post-stimulus events.
		}
	\label{fig:phase-reset-exp}
\end{figure}

An isolated oscillator (closed system) exhibits a precise rhythm, that is, a periodic behavior, and the period $T$ of the rhythm is assumed constant (see \myfigurename~\ref{fig:phase-reset-exp}a). 
To facilitate the comparison of rhythms with different periods (for example due to the variability in experimental preparation), it is convenient to define the notion of phase. In the absence of perturbations, the phase is a normalized time evolving on the unit circle.
Associating the onset of the observable event with phase $0$ (or $2\pi$), the phase variable $\theta(t)$ at time $t$ corresponds to the fraction of a period elapsed since the last occurrence of the observable event. It evolves linearly in time, that is, $\theta(t) \eqdef \omega \, (t - \hat{t}_i) \pmod{2\pi}$, where $\omega \eqdef 2\pi/T$ is the angular frequency of the oscillator and $\hat{t}_i$ is the time of the last observable event.

Following a phase-resetting stimulus at time $(t_s - \hat{t}_0)$ after one observable event (open system), the next event times $\hat{t}_i$, for $i\in\Nbb_{>0}$, are altered. For simplicity, it is assumed that the original rhythm is restored immediately after the first post-stimulus event, meaning that observable events repeat with the original period $T$ (see \myfigurename~\ref{fig:phase-reset-exp}b). 
The duration $\hat{T}\eqdef\hat{t}_1 - \hat{t}_0$ denotes the time interval from the event immediately before the stimulus to the next event after stimulation.
Once again, it is convenient to normalize each quantity in order to facilitate comparison between different experimental preparations. Multiplying by $\omega = 2\pi/T$ leads to $\theta \eqdef \omega\,(t_s - \hat{t}_0)$ and $\hat{\tau} \eqdef \omega\,\hat{T} = \omega\,(\hat{t}_1 - \hat{t}_0)$.

\begin{figure}[p]
	\centering
	\includegraphics{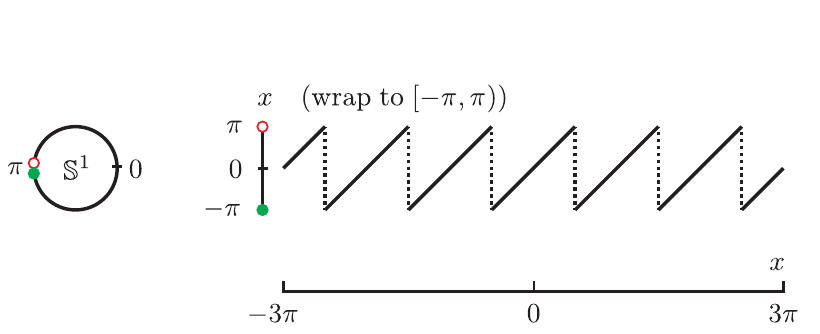}
	\caption{%
		Graphical representation of the wrap-to-$[-\pi,\pi)$ operation. 
		Given a real number $x$ in radians, $x \pwrap{[-\pi,\pi)} \equiv \left[x + \pi \pmod{2\pi}\right] - \pi$ wraps $x$ to the interval $[-\pi,\pi)$. 
		It adds or subtracts an integer multiple of $2\pi$ such that the result belongs to $[-\pi,\pi)$.
		(A solid dot indicates that the endpoint is included in the set; whereas an open dots indicates that the endpoint is excluded from the set.)
		}
	\label{fig:mod-wrap}
\end{figure}

The effect of a stimulation is to produce a phase shift~$\Delta \theta$ between the perturbed oscillator and the unperturbed oscillator. The phase shift $\Delta \theta$ is
\begin{equation*}
	\Delta \theta \eqdef 2\pi - \hat{\tau} \pwrap{[-\pi,\pi)},
\end{equation*}
where the operation $x \pwrap{[-\pi,\pi)} = [x + \pi \pmod{2\pi}] - \pi$ wraps $x$ to the interval $[-\pi,\pi)$ (see~\myfigurename~\ref{fig:mod-wrap}). 
Given a phase-resetting input $u(\cdot)$, the dependence of the phase shift $\Delta \theta$ on the (old) phase $\theta$ at which the stimulus was delivered is commonly called the phase response curve. It is denoted by $\PRC(\theta;u(\cdot))$, in order to stress that it is a function of the phase but that it also depends on the input $u(\cdot)$.

An alternative representation emphasizes the new phase~$\theta^+$ instead of the phase difference. Just before the stimulus, the oscillator had reached old phase $\theta$; just after, it appears to resume from the new phase $\theta^+$.  The new phase $\theta^+$ is 
\begin{equation*}
	\theta^+  \eqdef 2\pi - (\hat{\tau} - \theta) \pmod{2\pi}.
\end{equation*}
Given a phase-resetting input $u(\cdot)$, the dependence of the new phase $\theta^+$ on the (old) phase $\theta$ at which the stimulus was delivered is called the phase transition curve. It is denoted by $\PTC(\theta;u(\cdot))$.

Under the approximation that the initial rhythm is recovered immediately after the perturbation, the phase shift computed from the first post-stimulus event is identical to the asymptotic phase shift computed long after the perturbation.
This assumption neglects the transient change in the rhythm until a new steady-state is reached. 
To model the transient, the normalized time from the event before the stimulus to the $i$th event is denoted by $\hat{\tau}_i \eqdef \omega\,(\hat{t}_i - \hat{t}_0)$, leading to the phase shift $\Delta \theta_i \eqdef 2\pi - \hat{\tau}_i \pwrap{[-\pi,\pi)}$ and the new phase $\theta^+_i \eqdef 2\pi - (\hat{\tau}_i - \theta) \pmod{2\pi}$. 
If the oscillating phenomenon is time-invariant, a new steady-state behavior is expected asymptotically, such that $\lim_{i\rightarrow\infty} (\hat{\tau}_{i+1} - \hat{\tau}_{i}) = 2\pi$, $\lim_{i\rightarrow\infty}\Delta\theta_i \reveqdef \Delta\theta$, and $\lim_{i\rightarrow\infty}\theta^+_i \reveqdef \theta^+$. 

\clearpage
\section{Phase Response Curves from State-Space Models} \label{sec:prc_control}

This section reviews the mathematical characterization of phase response curves for oscillators described by time-invariant state-space models.

\subsection{State-Space Models of Oscillators}
Limit cycle oscillations appear in the context of nonlinear time-invariant state-space models
\begin{subequations} \label{eq:nlsys}
	\begin{align}
		\dot{x} & = \fvec(x,u), \\
			 y  & = \hvec(x),
	\end{align}	
\end{subequations}
where the states $x(t)$ evolve on some subset $\mathcal{X}\subseteq\Rbb^n$, and the input and output values $u(t)$ and $y(t)$ belong to subsets $\mathcal{U}\subseteq\Rbb$ and $\mathcal{Y}\subseteq\Rbb$, respectively. 
The vector field $\fvec:\mathcal{X}\times\mathcal{U} \rightarrow \Rbb^n$ and the measurement map $\hvec:\mathcal{X} \rightarrow \mathcal{Y}$ support all the usual smoothness conditions that are necessary for existence and uniqueness of solutions. 
An input is a signal $u:[0,\infty)\rightarrow\mathcal{U}$ that is locally essentially compact (meaning that images of restrictions to finite intervals are compact).
The solution at time $t$ to the initial value problem $\dot{x}  = \fvec(x,u)$ from the initial condition $x_0\in\mathcal{X}$ at time $0$ is denoted by $\flow(t,x_0,u(\cdot))$ (with $\flow(0,x_0,u(\cdot)) = x_0$). 
For convenience, single-input and single-output systems are considered. 
All developments generalize to multiple-input and multiple-output systems.

The state-space model~\eqref{eq:nlsys} is called an oscillator if the zero-input system $\dot{x}=\fvec(x,0)$ admits an exponentially stable limit cycle, that is, a periodic orbit $\perorb\subseteq\mathcal{X}$ with period $T$ that attracts nearby solutions at an exponential rate~\cite{Farkas:1994uq}.
Picking an initial condition $\persol_0\in\perorb$, the periodic orbit $\perorb$ is described by the locus of the (nonconstant) $T$-periodic solution $\flow(\cdot,\persol_0,\zeroinput)$, that is,
\begin{equation*}
	\perorb \eqdef \left\{ x \in \mathcal{X} : x = \flow(t,\persol_0,\zeroinput), t \in [0,T) \right\},
\end{equation*}
where the period $T>0$ is the smallest positive constant such that $\flow(t,\persol_0,\zeroinput) = \flow(t+T,\persol_0,\zeroinput)$ for all $t\geq 0$ and  $\zeroinput$ is the input signal identically equal to $0$ for all times. The periodic orbit is an invariant set.

Because of the periodic nature of the steady-state behavior, it is appealing to study the oscillator dynamics directly on the unit circle $\Sbb^1$. The key ingredient of this phase reduction is the phase map concept. 
A phase map $\Theta:\mathcal{B}(\gamma)\subseteq\mathcal{X}\rightarrow\Sbb^1$ is a mapping that associates with every point in the basin of attraction $\mathcal{B}(\gamma)\subseteq\mathcal{X}$ a phase on the unit circle $\Sbb^1$.  Away from a finite number of isolated points (called singular points), the phase map $\Theta$ is a continuous map. 
The phase variable $\theta(t)$ is the image of the flow through the phase map, that is, $\theta(t)\eqdef\Theta(\flow(t,x_0,u(\cdot)))$. 
By the definition of the phase map, the phase dynamics reduce to $\dot{\theta} = \omega$ for the input $\zeroinput$.
For nonzero inputs, the phase dynamics are often hard to derive. See ``Phase Maps'' for details.

For convenience, the periodic orbit $\perorb$ is parameterized by the map $\persol:\Sbb^1\rightarrow\perorb$ that associates with each phase $\theta$ on the unit circle a point $\flow(\theta/\omega,\persol_0,\zeroinput)\reveqdef\persol(\theta)$ on the periodic orbit.

\subsection{Response to Phase-Resetting Inputs}

If a solution of \eqref{eq:nlsys} asymptotically converges to the periodic orbit, the corresponding input~$u(\cdot)$ is said to be phase-resetting. If an input is phase-resetting for an initial condition~$x_0$, then there exists a phase shift $\theta^+\in\Sbb^1$ that satisfies
\begin{equation*} \label{eq:reset_phase}
	\lim_{t\rightarrow\infty} \left\| \flow(t,x_0,u(\cdot)) - \flow(t , \flow(\theta^+/\omega,\persol_0,\zeroinput),\zeroinput) \right\|_2 = 0 .
\end{equation*}

\begin{defn} \label{def:prc}
	Given a phase-resetting input $u(\cdot)$, the (finite) phase response curve is the map $\PRC(\cdot;u(\cdot)):\Sbb^1 \rightarrow [-\pi,\pi)$ that associates with each phase $\theta$ a phase shift $\Delta\theta = \PRC(\theta;u(\cdot))$, defined as
	\begin{align*}
		\PRC(\theta;u(\cdot)) 
				& = \lim_{t \rightarrow +\infty} [\Theta(\flow(t,\persol(\theta),u(\cdot))) - (\omega\,t + \theta)] \pwrap{[-\pi,\pi)} . \\
	\end{align*}
	Similarly, the phase transition curve is the map  $\PTC(\cdot;u(\cdot)):\Sbb^1 \rightarrow \Sbb^1$ that associates with each phase $\theta$ the new phase $\theta^+ = \PTC(\theta;u(\cdot))$, defined as
	\begin{align*}
	\PTC(\theta;u(\cdot)) & = \lim_{t \rightarrow +\infty} [\Theta(\flow(t,\persol(\theta),u(\cdot))) - \omega\,t] \pmod{2\pi}. \\
	\end{align*}
\end{defn}

A mathematically more abstract---yet useful---tool is the infinitesimal phase response curve. It captures the same information as the finite phase response curve in the limit of Dirac delta input with infinitesimal amplitude (that is, $u(\cdot) = \alpha\,\delta(\cdot)$ with $\alpha \rightarrow 0$). 
\begin{defn}
	The infinitesimal phase response curve is the map $\iPRCu:\mathbb{S}^1\rightarrow\mathbb{R}$, defined as the directional derivative 
	\begin{equation*} \label{eq:q_u}
		\iPRCu(\theta) \eqdef \dirder{\Theta}{\persol(\theta)}{\frac{\partial \fvec}{\partial u}(\persol(\theta),0)} ,
	\end{equation*}
	where
	\begin{equation*}
		\dirder{\Theta}{x}{\eta} \eqdef \lim_{h\rightarrow0}\frac{\Theta(x+h\,\eta) - \Theta(x)}{h}.
	\end{equation*}
\end{defn}
The directional derivative can be computed as the inner product in $\Rbb^n$
	\begin{equation} \label{eq:dirder}
		\iPRCu(\theta) = \dirder{\Theta}{\persol(\theta)}{\frac{\partial \fvec}{\partial u}(\persol(\theta),0)} = \left\langle \Euclgrad{\Theta}{\persol(\theta)}{x} , \frac{\partial \fvec}{\partial u}(\persol(\theta),0) \right\rangle ,
	\end{equation}
	where $\Euclgrad{\Theta}{\persol(\theta)}{x} \reveqdef p(\theta)$ is the gradient of the asymptotic phase map $\Theta$ at the point $\persol(\theta)$. 
	
The main benefit of an infinitesimal characterization of phase response curves is that the concept is independent of the input signal. Limitations of the infinitesimal approach have been well identified since the early days of phase resetting studies~\cite{Winfree:1980ue} and strongly depend on the application context. For instance, infinitesimal phase response curves have proven very useful in the study of circadian rhythms~\cite{Pfeuty:2011em}, but come with severe limitations in the context of neurodynamics, as recently studied in~\cite{Oprisan:2003ki,Achuthan:2009tu,Wang:2013ab}.
	
\begin{rem}
	By definition, the finite phase response curve for an impulse input is well approximated by the infinitesimal phase response curve, that is, $\PRC(\cdot;\alpha\,\delta(\cdot)) = \alpha \, q(\cdot) + \mathcal{O}(\alpha^2)$. See ``From Infinitesimal to Finite Phase Response Curves'' for details.
\end{rem}
	
\subsection{Phase Models as Reduced Models of Oscillators}

Phase response curves are the basis for the reduction of $n$-dimensional state-space models of oscillators to one-dimensional phase models. 
Phase models are the main representation of oscillators for networks. However, the focus of this article is on single oscillator models.
For a comprehensive treatment of phase models, the reader is referred to the vast literature on the subject (see pioneering papers \cite{Winfree:1967vf,Kuramoto:1975ki,Ermentrout:1984jb,Mirollo:1990ft}, review articles~\cite{Mauroy:2012vi,Dorfler:2012wv,Dorfler:2013fk}, and books \cite{Winfree:1980ue,Glass:1988ub,Hoppensteadt:1997tp,Izhikevich:2007vr,Schultheiss:2012bz}).


Below, two popular phase models are reviewed. They are obtained through phase reduction methods in the case of weak input and impulse train input, respectively.

Under the simplifying assumption of weak input, that is,
\begin{equation*}
	|u(t)| \ll 1, \quad \text{for all $t\geq0$},
\end{equation*}
any solution $\flow(t,x_0,u(\cdot))$ of the oscillator model that starts in the neighborhood of the hyperbolic stable periodic orbit~$\gamma$ stays in its neighborhood. The $n$-dimensional state-space model \eqref{eq:nlsys} can thus be approximated by a one-dimensional continuous-time phase model (see~\cite{Kuramoto:1984wo,Kuramoto:1997kd,Hoppensteadt:1997tp,Brown:2004iy,Izhikevich:2007vr,Mauroy:2012vi})
\begin{subequations}\label{eq:weak_coupling}
	\begin{align} 
		\dot{\theta} & = \omega + \iPRCu(\theta) \, u, \\
		y            & = \tilde{\hvec}(\theta),
	\end{align}	
\end{subequations}
where the phase variable $\theta$ evolves on the unit circle $\mathbb{S}^1$. The phase model is fully characterized by the angular frequency $\omega>0$, the infinitesimal phase response curve $\iPRCu:\mathbb{S}^1 \rightarrow \mathbb{R}$, and the measurement map $\tilde{h} :\Sbb^1\rightarrow\mathcal{Y}$, which is defined as $\tilde{h}(\theta) = h(\persol(\theta))$.

An alternative simplification is when the input is a train of resetting impulses, that is, 
\begin{equation*}
	u(t) = \alpha \, \sum_{k=0}^{\infty} \delta(t-t_{k}), \quad \text{with $t_k \geq 0$},
\end{equation*}
where it is assumed that the time interval between successive impulses is sufficient for convergence to the periodic orbit between each of them. Under this assumption, any solution $\flow(t,x_0,u(\cdot))$ of the oscillator model that starts from the periodic orbit $\gamma$ leaves the periodic orbit under the effect of one impulse from the train and then converges back toward the periodic orbit.
Assuming that the steady-state of the periodic orbit is recovered between any two successive impulses, the $n$-dimensional state-space model \eqref{eq:nlsys} can be approximated by a one-dimensional hybrid phase model (see~\cite{Glass:1988ub,Izhikevich:2007vr,Mauroy:2012vi}) with
\begin{subequations} \label{equa_hybrid}
	\begin{enumerate}
		\item the (constant-time) flow rule 
		\begin{align}
			\dot{\theta} & = \omega, && \text{for all $t \neq t_k$}, \\
		\intertext{\item the (discrete-time) jump rule}
			\theta^+  & = \theta + \PRC(\theta;\alpha\,\delta(\cdot)), && \text{for all $t = t_k$}, \\
		\intertext{\item the measurement map}	
			y          & = \tilde{\hvec}(\theta), && \text{for all $t$} ,
		\end{align}
	\end{enumerate}	
\end{subequations}
where the phase variable $\theta$ evolves on the unit circle $\mathbb{S}^1$. The phase model is fully characterized by the angular frequency $\omega>0$, the phase response curve $\PRC(\cdot;\alpha\,\delta(\cdot)):\mathbb{S}^1\rightarrow[-\pi,\pi)$, and the measurement map $\tilde{h} :\Sbb^1\rightarrow\mathcal{Y}$.

It should be emphasized that the assumption of ``weak inputs'' or ``trains of resetting impulses'' is relative to the attractivity of the periodic orbit. Strongly attractive periodic orbits allow for larger inputs to meet the simplifying assumption. The use of phase models is for instance popular in the study of oscillator networks under the assumption that the coupling strength is weak with respect to the attractivity of each oscillator \cite{Mauroy:2012vi,Dorfler:2012wv,Dorfler:2013fk}.

Both reduced oscillator representations $\{\omega,\iPRCu(\cdot),\tilde{\hvec}(\cdot)\}$ and $\{\omega,\PRC(\cdot;\alpha\,\delta(\cdot)),\tilde{\hvec}(\cdot)\}$ have characteristics similar to the static gain in the transfer-function representation of linear time-invariant systems. Both representations capture asymptotic properties of the impulse response. They are external input--output representations of the oscillators, independent of the complexity of the internal state-space representation of the oscillators. Moreover, information on such characteristics is available experimentally.

\subsection{Computations of Phase Response Curves}

A brief review of numerical methods to compute periodic orbits and phase response curves in state-space models is useful before introducing the numerics of sensitivity analysis.

\subsubsection{Periodic Orbit}

The $2\pi$-periodic steady-state solution $\persol(\cdot)$ and the angular frequency $\omega$ are calculated by solving the boundary value problem (see  \cite{Ascher:1988ty} and \cite{Seydel:2010tw})
\begin{subequations} \label{eq:bvp_x}
	\begin{align}
		\frac{d\persol}{d\theta}(\theta) - \frac{1}{\omega} \, \fvec(\persol(\theta),0) & = 0, \label{eq:bvp_x_dx} \\
		\persol(2\pi) - \persol(0) & = 0, \label{eq:bvp_x_per}\\
		\PLC(\persol(0))  & = 0. \label{eq:bvp_x_phase}
	\end{align}	
\end{subequations}
The boundary conditions are given by the periodicity condition \eqref{eq:bvp_x_per}, which that ensures the periodicity of the map $\persol(\cdot)$, and the phase condition \eqref{eq:bvp_x_phase}, which anchors a reference position $\persol(0) = \persol_0$ along the periodic orbit. The phase condition $\PLC:\mathcal{X}\rightarrow\Rbb$ is chosen such that it vanishes at an isolated point $\persol_0$ on the periodic orbit $\perorb$ (see~\cite{Seydel:2010tw} for details).
Numerical algorithms to solve this boundary value problem are reviewed in ``Numerical Tools''.

\subsubsection{Infinitesimal Phase Response Curve}

The infinitesimal phase response curve $\iPRCu(\cdot)$ is calculated by applying \eqref{eq:dirder} that involves computing the gradient of the asymptotic phase map evaluated along the periodic orbit, that is, the function~$p(\cdot)$.

The gradient of the asymptotic phase map evaluated along the periodic orbit $\iPRCx(\cdot)$ is calculated by solving the boundary value problem (see \cite{Malkin:1949to,Malkin:1956wq,Neu:1979cj,Ermentrout:1984jb,Kuramoto:1984wo,Govaerts:2006bv})
\begin{subequations} \label{eq:bvp_q}
	\begin{align}
		\frac{d\iPRCx}{d\theta}(\theta) + \frac{1}{\omega} \, \trans{\frac{\partial \fvec}{\partial x}(\persol(\theta),0)} \, \iPRCx(\theta) & = 0, \label{eq:bvp_q_a} \\
		\iPRCx(2\pi) - \iPRCx(0) & = 0, \label{eq:bvp_q_b} \\
		\langle \iPRCx(\theta) , \fvec(\persol(\theta),0) \rangle - \omega & = 0, \label{eq:bvp_q_c}
	\end{align}	
\end{subequations}
where the notation $\trans{A}$ stands for the transpose of the matrix $A$. The boundary condition \eqref{eq:bvp_q_b} imposes the periodicity of~$\iPRCx(\cdot)$ and the normalization condition \eqref{eq:bvp_q_c}~ensures a linear increase at rate $\omega$ of the phase variable~$\theta$ along zero-input trajectories. 
This method is often called the adjoint method.
Numerical methods to solve this boundary value problem as a by-product of the periodic orbit computation are presented in ``Numerical Tools''.

\subsubsection{Finite Phase Response Curve}

As an alternative to the infinitesimal phase response curve, direct methods compute numerically the phase response curve of an oscillator state-space model as a direct application of Definition~\ref{def:prc} (see for example \cite{Winfree:1974vm,Winfree:1980ue,Glass:1988ub,Kuramoto:1997kd,Brown:2004iy,Taylor:2008co}). 

For each point $\PRC(\theta_i;u(\cdot))$, with  $\theta_i\in\Sbb^1$, of the finite phase response curve, a perturbed trajectory $\flow(t,\persol(\theta_i),u(\cdot))$ is computed by solving the initial value problem~\eqref{eq:nlsys} from $\persol(\theta_i)$ up to its convergence back in a neighborhood of the periodic orbit, that is, up to time $t_*$ such that $\dist(\flow(t_*,\persol(\theta_i),u(\cdot)),\gamma) < \epsilon$, where $\epsilon$ is a chosen  error tolerance. The phase $\theta_* = \Theta(\flow(t_*,\persol(\theta_i),u(\cdot)))$  is estimated as
\begin{equation*}
	\theta_* = \argmin_{\theta\in\Sbb^1} \left\|\flow(t_*,\persol(\theta_i),u(\cdot)) - \persol(\theta)\right\|_2.
\end{equation*}
Then, the asymptotic phase shift is measured by direct comparison with the phase $\omega\,t_*+\theta_i$ of an unperturbed trajectory at time $t_*$, that is, 
\begin{equation*}
	\PRC(\theta_i;u(\cdot)) = \theta_* - (\omega\,t_*+\theta_i).
\end{equation*}

An advantage of the direct method over the infinitesimal method is that it applies to arbitrary phase-resetting inputs.
It only requires an efficient time integrator. However, it is highly expensive from a computational point of view: for each phase-resetting input, each point of the corresponding phase response curve requires the time simulation of the $n$-dimensional state-space model, up to the asymptotic convergence of the perturbed trajectory towards the periodic orbit.

\clearpage
\section{Metrics in the Space of Phase Response Curves} \label{sec:metrics}

To answer system-theoretic questions in the space of phase response curves, it is useful to endow this space with the differential structure of a Riemannian manifold. The differential structure provides a notion of local sensitivity in the tangent space. The Riemannian structure is convenient for recasting analysis problems in an optimization framework because it provides, for instance, a notion of steepest descent. The Riemannian structure also provides a norm in the tangent space and a (geodesic) distance between phase response curves. See ``Basic Concepts of Differential Geometry on Manifolds'' for a short introduction to these concepts.

Because phase response curves are signals defined on the unit circle and take values on the real line, the most obvious Riemannian structure is provided by the infinite-dimensional Hilbert space of square-integrable signals
\begin{equation*}
	\Hilbert^{0}  \eqdef \{q:q(\cdot)\in\Lp{2}{\Sbb^1}{\Rbb}\},
\end{equation*}
where $\Lp{2}{\Sbb^1}{\Rbb} = \{ q:\Sbb^1\rightarrow\Rbb : (\int_{0}^{2\pi} |q(\theta)|^2 \, d\theta)^{\frac{1}{2}} < \infty \}$,
endowed with the standard inner product
\begin{equation} \label{eq:H0-inner-product}
	\langle \xi(\cdot) , \zeta(\cdot) \rangle \eqdef \int_{0}^{2\pi} \xi(\theta) \, \cconj{\zeta(\theta)} \, d\theta
\end{equation}
and the associated norm
\begin{equation} \label{eq:H0-norm}
	\norm{\xi(\cdot)}_2 \eqdef \sqrt{\langle \xi(\cdot) , \xi(\cdot) \rangle}.
\end{equation}

For technical reasons detailed later, the first derivative of considered signals is also assumed to be square-integrable. It thus restricts the signal space to 
\begin{equation*}
	\Hilbert^{1}  \eqdef \left\{ q : q(\cdot) \in \Lp{2}{\Sbb^1}{\Rbb}, q'(\cdot) \in \Lp{2}{\Sbb^1}{\Rbb} \right\},
\end{equation*}
where $q'$ denotes the derivative, with respect to the phase $\theta$, of the signal $q$.
The space $\Hilbert^{1}$ is a linear subspace of $\Hilbert^{0}$ and it inherits its inner product \eqref{eq:H0-inner-product} and its norm \eqref{eq:H0-norm}.

The linear space structure $\Hilbert^{1}$ is convenient for calculations but it fails to capture natural equivalence properties between phase response curves. In many applications, it is not meaningful to distinguish among phase response curves that are related by a scaling factor and/or a  phase shift.

\paragraph*{Scaling equivalence}

The actual magnitude of the input signal acting on the system is not always known exactly. This uncertainty about the input magnitude induces an (inversely proportional) uncertainty about the phase response magnitude. 
Indeed, the phase model \eqref{eq:weak_coupling} is equivalent to 
\begin{align*}
	\dot{\theta} & = \omega + \left(q(\theta)\,\alpha\right) \, \left(\frac{1}{\alpha}\,u\right), \\
	y & = \tilde{\hvec}(\theta),
\end{align*}
for any scaling factor $\alpha>0$. A scaling of the input magnitude can be counterbalanced by an inverse scaling of the phase response curve.
In these cases, a phase response curve $q$ is considered as the representation of an equivalence class $\sim$ characterized by
\begin{equation} \label{eq:scaling-equiv}
	q_1 \sim q_2 \Leftrightarrow \text{there exists $\alpha>0:q_2(\cdot) = q_1(\cdot) \, \alpha$}  .
\end{equation}

For example, in circadian rhythms, the stimulus could be a pulse of light, the effect of drugs, or the intake of food.  Pulses are modeled by scaling the intensity of a parameter but the absolute variation of this parameter is not known and is empirically fitted to experimental data. The scaling equivalence is meaningful in such situations. 
On the other hand, in neurodynamics, the stimulus could be a post-synaptic current of constant magnitude. In this latter case, the scaling equivalence is less appropriate.

\paragraph*{Phase-shifting equivalence}

The choice of a reference position (associated with the initial phase) along the periodic orbit is often arbitrary.
In these cases, a phase response curve $q$ is considered as representative of an equivalence class $\sim$ characterized by
\begin{equation} \label{eq:shifting-equiv}
	q_1 \sim q_2 \Leftrightarrow \text{there exists $\sigma\in\Sbb^1:q_2(\cdot) = q_1(\cdot + \sigma)$},
\end{equation}
where $\sigma$ denotes any phase shift. 

For example, in circadian rhythms, experimental data are often collected by observing the locomotor activity of the animal. The timing of this locomotor activity is not easily linked to the time evolution of molecular concentrations. In this case, the phase shifting equivalence is meaningful. 
On the other hand, in neurons, the observable events are the action potentials measured as rapid changes in membrane potentials. If the membrane potential is a state variable of the model, there is no timing ambiguity. In this latter case, the phase-shifting equivalence is not appropriate.

The equivalence relations \eqref{eq:scaling-equiv} and \eqref{eq:shifting-equiv} lead to the abstract---yet useful---concept of quotient space. Each point of a quotient space is defined as an equivalence class of signals. Since these equivalence classes are abstract objects, they cannot be used explicitly in numerical computations. Algorithms on quotient space work instead with representatives (in the total space) of these equivalence classes.

Combining (or not) equivalence properties \eqref{eq:scaling-equiv} and \eqref{eq:shifting-equiv} ends up with four infinite-dimensional spaces: one Hilbert space and three quotient spaces, respectively, denoted by $\mathcal{Q}_{\text{A}}$, $\mathcal{Q}_{\text{B}}$, $\mathcal{Q}_{\text{C}}$, and $\mathcal{Q}_{\text{D}}$ (see \mytablename~\ref{tab:metric}). 
In the next four subsections, each space is endowed with an appropriate Riemannian metric and an expression of tangent vectors, needed for the sensitivity analysis in subsequent sections, is provided.

Below, the symbol $q$ denotes an element of the considered space. It can be a signal (a finite or infinitesimal phase response curve) or an equivalence class of these signals. In the later case, a signal is denoted by $\bar{q}$.

\subsection{Metric on Hilbert Space $\Hilbert^{1}$}

The simplest space structure is the Hilbert space $\mathcal{Q}_\text{A} \eqdef \Hilbert^{1}$.
The (flat) Riemannian metric on $\mathcal{Q}_\text{A}$ is the inner product
\begin{equation*}
	g_q(\xi_q,\zeta_q) \eqdef \langle \xi_q,\zeta_q \rangle
\end{equation*}
with (Euclidean) induced norm
\begin{equation*}
	\norm{\xi_q}_q \eqdef \sqrt{g_q(\xi_q,\xi_q)} = \sqrt{\langle\xi_q,\xi_q\rangle} = \norm{\xi_q}_2.
\end{equation*}

Because the space $\mathcal{Q}_{\text{A}}$ is a linear space structure, the shortest path between two elements $q_1$ and $q_2$ on $\mathcal{Q}_{\text{A}}$ is the straight line joining these elements. The natural (geodesic) distance between two points $q_1$ and $q_2$ on $\mathcal{Q}_{\text{A}}$ is then given by
\begin{equation*}
	\dist(q_1,q_2) \eqdef \norm{ q_1-q_2 }_{2}.
\end{equation*}

\subsection{Metric on the Quotient Space $\Hilbert^{1}/\Rbb_{>0}$}

The space capturing the scaling equivalence \eqref{eq:scaling-equiv} is the quotient space $\mathcal{Q}_\text{B}\eqdef\Hilbert^{1}/\Rbb_{>0}$. Each element $q$ in $\mathcal{Q}_\text{B}$ represents an equivalence class
\begin{equation*}
	q = [\overline{q}] \eqdef \{ \overline{q} \, \alpha : \alpha > 0 \}.
\end{equation*}
These equivalence classes are rays (starting at $0$) in the total space $\overline{\mathcal{Q}}_\text{B}\eqdef\Hilbert^{1}$.

The normalized metric on $\overline{\mathcal{Q}}_\text{B}$,
\begin{equation} \label{eq:metric_B}
	\overline{g}_{\overline{q}}(\overline{\xi}_{\overline{q}},\overline{\zeta}_{\overline{q}}) \eqdef \frac{\langle \overline{\xi}_{\overline{q}},\overline{\zeta}_{\overline{q}} \rangle}{\langle \overline{q},\overline{q} \rangle},
\end{equation} 
is invariant by scaling. As a consequence, it induces a Riemannian metric $g_{q}(\xi_{q},\zeta_{q}) \eqdef \overline{g}_{\overline{q}}(\overline{\xi}_{\overline{q}},\overline{\zeta}_{\overline{q}})$ on $\mathcal{Q}_\text{B}$.
The norm in the tangent space $\tgspace{q}{\mathcal{Q}_\text{B}}$ at $q$ is 
\begin{equation} \label{eq:normB}
	\norm{\xi_q}_q \eqdef \sqrt{g_q(\xi_q,\xi_q)} = \frac{\norm{\overline{\xi}_{\overline{q}}}_2}{\norm{\overline{q}}_2}.
\end{equation}

A signal representation of a tangent vector at $q\in\mathcal{Q}_\text{B}$ relies on the decomposition of the tangent space $\tgspace{\overline{q}}{\overline{\mathcal{Q}}_\text{B}}$ into its vertical and horizontal subspaces. The vertical subspace $\mathcal{V}_{\overline{q}}$ is the subspace of $\tgspace{\overline{q}}{\overline{\mathcal{Q}}_\text{B}}$ that is tangent to the equivalence class $[\overline{q}]$, that is,
\begin{equation*}
	\mathcal{V}_{\overline{q}} = \{\overline{q} \, \beta:\beta\in\Rbb\}.
\end{equation*}
The horizontal space $\mathcal{H}_{\overline{q}}$ is chosen as the orthogonal complement of $\mathcal{V}_{\overline{q}}$ for the metric $\overline{g}_{\overline{q}}(\cdot,\cdot)$, that is, 
\begin{equation*}
	\mathcal{H}_{\overline{q}} = \{\eta \in \tgspace{\overline{q}}{\overline{\mathcal{Q}}_\text{B}}:\overline{g}_{\overline{q}}(\eta,\overline{q} \, \beta) = 0\}.
\end{equation*}
The orthogonal projection $P_{\overline{q}}^h \, \eta$ of a vector $\eta\in\tgspace{\overline{q}}{\overline{\mathcal{Q}}_\text{B}}$ onto the horizontal space $\mathcal{H}_{\overline{q}}$ is
\begin{equation*}
	P_{\overline{q}}^h \, \eta \eqdef \eta - \frac{\overline{g}_{\overline{q}}(\eta,\overline{q}\,\beta)}{\overline{g}_{\overline{q}}(\overline{q}\,\beta,\overline{q}\,\beta)} \, \overline{q} \, \beta = \eta - \frac{\langle\eta,\overline{q}\rangle}{\langle\overline{q},\overline{q}\rangle} \, \overline{q} .
\end{equation*}

The distance between two points $q_1$ and $q_2$ on $\mathcal{Q}_\text{B}$ is defined as
\begin{equation*}
	\dist(q_1,q_2) \eqdef \cos^{-1}\left( \frac{\langle \overline{q}_1, \overline{q}_2 \rangle}{\norm{\overline{q}_1}_2 \, \norm{\overline{q}_2}_2} \right)
\end{equation*}
(see~\cite{Goh:2008tc} for metrics on the unit sphere).

\subsection{Metric on the Quotient Space $\Hilbert^{1}/\Shift(\Sbb^1)$}

The space capturing the phase-shifting equivalence \eqref{eq:shifting-equiv} is the quotient space $\mathcal{Q}_\text{C}\eqdef\Hilbert^{1}/\Shift(\Sbb^1)$.
Each element $q$ in $\mathcal{Q}_\text{C}$ represents an equivalence class
\begin{equation*}
	q = [\overline{q}] = \{ \overline{q}(\cdot+\sigma) : \sigma \in \Sbb^1 \}.
\end{equation*}
These equivalence classes are closed one-dimensional curves (due to the periodicity of the shift) on the infinite-dimensional hypersphere of radius $\norm{\overline{q}}_2$ in the total space $\overline{\mathcal{Q}}_\text{C}\eqdef\Hilbert^{1}$.

The (flat) metric on $\overline{\mathcal{Q}}_\text{C}$
\begin{equation*}
	\overline{g}_{\overline{q}}(\overline{\xi}_{\overline{q}},\overline{\zeta}_{\overline{q}}) \eqdef \langle \overline{\xi}_{\overline{q}},\overline{\zeta}_{\overline{q}} \rangle,
\end{equation*}
is invariant by phase shifting along the equivalence classes. As a consequence, it induces a Riemannian metric $g_{q}(\xi_{q},\zeta_{q}) \eqdef \overline{g}_{\overline{q}}(\overline{\xi}_{\overline{q}},\overline{\zeta}_{\overline{q}})$ on $\mathcal{Q}_\text{C}$.
The norm in the tangent space $\tgspace{q}{\mathcal{Q}_\text{C}}$ at $q$ is 
\begin{equation*}
	\norm{\xi_q}_q \eqdef \sqrt{g_q(\xi_q,\xi_q)} = \norm{\overline{\xi}_{\overline{q}}}_2.
\end{equation*}

The vertical space $\mathcal{V}_{\overline{q}}$ is the subspace of $\tgspace{\overline{q}}{\overline{\mathcal{Q}}}_\text{C}$ that is tangent to the equivalence class $[\overline{q}]$, that is,
\begin{equation*}
	\mathcal{V}_{\overline{q}} = \{\overline{q}' \,\beta:\beta\in\Rbb\},
\end{equation*}
where $\overline{q}'$ has to belong to $\Lp{2}{\Sbb^1}{\Rbb}$ to ensure the regularity of $\mathcal{V}_{\overline{q}}$.
The horizontal space $\mathcal{H}_{\overline{q}}$ is chosen as the orthogonal complement of $\mathcal{V}_{\overline{q}}$ for the metric $\overline{g}_{\overline{q}}(\cdot,\cdot)$, that is,
\begin{equation*}
	\mathcal{H}_{\overline{q}} = \{\eta \in \tgspace{\overline{q}}{\overline{\mathcal{Q}}_{\text{C}}}:\overline{g}_{\overline{q}}(\eta , \overline{q}' \, \beta ) = 0\}.
\end{equation*}
The orthogonal projection $P_{\overline{q}}^h \,\eta$ of a vector $\eta\in\tgspace{\overline{q}}{\overline{\mathcal{Q}}_\text{C}}$ onto the horizontal space $\mathcal{H}_{\overline{q}}$ is
\begin{equation*}
	P_{\overline{q}}^h \, \eta \eqdef \eta - \frac{\overline{g}_{\overline{q}}(\eta,\overline{q}'\,\beta)}{\overline{g}_{\overline{q}}(\overline{q}'\,\beta,\overline{q}'\,\beta)}\, \overline{q}'\,\beta = \eta - \frac{\langle\eta,\overline{q}'\rangle}{\langle\overline{q}',\overline{q}'\rangle} \, \overline{q}' .
\end{equation*}

The distance between two points $q_1$ and $q_2$ on $\mathcal{Q}_\text{C}$ is defined as
\begin{equation*}
	\dist(q_1,q_2) \eqdef \min_{\sigma\in\Sbb^1} \norm{\overline{q}_1(\cdot)-\overline{q}_2(\cdot+\sigma)}_2  = \norm{\overline{q}_1(\cdot)-\overline{q}_2(\cdot+\sigma_*)}_2,
\end{equation*}
where $\sigma_*$ denotes the phase shift achieving this minimization. It corresponds to the phase shift  maximizing the circular cross-correlation
\begin{equation} \label{eq:opt-phase-shift}
	\sigma_* = \argmax_{\sigma \in \Sbb^1} \langle \overline{q}_1(\cdot),\overline{q}_2(\cdot+\sigma) \rangle.
\end{equation}

The global optimization problem \eqref{eq:opt-phase-shift} is solved in two steps.
The first step is the computation of the circular cross-correlation $\overline{c}:\Sbb^1\rightarrow\Rbb$ between the two periodic signals $\overline{q}_1$ and $\overline{q}_2$
	\begin{equation*}
		\overline{c}(\sigma) = \langle \overline{q}_1(\cdot),\overline{q}_2(\cdot+\sigma) \rangle.
	\end{equation*}
	By definition, the circular cross-correlation is also a periodic signal.
	An efficient computation of this circular cross-correlation is performed in the Fourier domain. The circular cross-correlation can be expressed as the circular convolution $\overline{c}(\sigma) = ( \cconj{\overline{q}_1(-\cdot)} \odot \overline{q}_2(\cdot) ) (\sigma)$. Exploiting the properties of Fourier coefficients and the convolution-multiplication duality property leads to
	\begin{equation*}
		\hat{\overline{c}}[k] = \cconj{\hat{\overline{q}}_1[k]} \, \hat{\overline{q}}_2[k]	,
	\end{equation*}
	where $\hat{x}[\cdot]$ denotes the discrete signal of Fourier coefficients for the periodic signal $x(\cdot)$ and $\cconj{x}$ denotes the complex conjugate of $x$.
The second step is the identification of the optimal phase-shift value $\sigma_* \in \Sbb^1$, which achieves the maximal value of the circular cross-correlation. This maximum is global and generically unique. Multiplicity of the optimum would mean that one of the signals has a period that is actually equal to $2\pi/k$ with $k\in\Nbb_{>0}$.

\subsection{Metric on the Quotient Space $\Hilbert^{1}/(\Rbb_{>0}\times\Shift(\Sbb^1))$}

The space capturing both scaling and phase-shifting equivalences \eqref{eq:scaling-equiv}--\eqref{eq:shifting-equiv} is the quotient space $\mathcal{Q}_\text{D}\eqdef\Hilbert^{1}/(\Rbb_{>0}\times\Shift(\Sbb^1))$. Each element $q$ in $\mathcal{Q}_\text{D}$ represents an equivalence class
\begin{equation*}
	q = [\overline{q}] = \{ \overline{q}(\cdot+\sigma) \, \alpha : \alpha>0, \sigma \in \Sbb^1 \}.
\end{equation*}
Based on the individual geometric interpretation of both equivalence properties,  these equivalence classes are infinite cones in the total space $\overline{\mathcal{Q}}_\text{D}\eqdef\Hilbert^1$, that is, the union of rays that start at $0$ and go through the closed one-dimensional curve of phase-shifted signals.

Because the metric \eqref{eq:metric_B} on $\overline{\mathcal{Q}}_\text{D}$ is invariant by scaling and phase shifting along the equivalence classes, it induces a Riemannian metric $g_{q}(\xi_{q},\zeta_{q}) \eqdef \overline{g}_{\overline{q}}(\overline{\xi}_{\overline{q}},\overline{\zeta}_{\overline{q}})$ on $\mathcal{Q}_\text{D}$.
The norm in the tangent space $\tgspace{q}{\mathcal{Q}_\text{D}}$ at $q$ is given by \eqref{eq:normB}.

The vertical space $\mathcal{V}_{\overline{q}}$ is the subspace of $\tgspace{\overline{q}}{\overline{\mathcal{Q}}}_\text{D}$ that is tangent to the equivalence class $[\overline{q}]$, that is,
\begin{equation*}
	\mathcal{V}_{\overline{q}} = \{\overline{q}\,\beta_1+\overline{q}'\,\beta_2:\beta_1,\beta_2\in\Rbb\}.
\end{equation*}
It is the direct sum of vertical spaces for equivalence properties individually. 
The horizontal space $\mathcal{H}_{\overline{q}}$ is chosen as the orthogonal complement of $\mathcal{V}_{\overline{q}}$ for the metric $g_q(\cdot,\cdot)$, that is,
\begin{equation*}
	\mathcal{H}_{\overline{q}} = \{\eta \in \tgspace{\overline{q}}{\overline{\mathcal{Q}}_{\text{D}}}: \overline{g}_{\overline{q}}(\eta,\overline{q}\,\beta_1+\overline{q}'\,\beta_2) = 0 \}.
\end{equation*}
The orthogonal projection $P_{\overline{q}}^h \, \eta$ of a vector $\eta\in\tgspace{\overline{q}}{\overline{\mathcal{Q}}_\text{D}}$ onto the horizontal space $\mathcal{H}_{\overline{q}}$ is
\begin{align*}
	P_{\overline{q}}^h \, \eta 
	& \eqdef  \eta - \frac{\overline{g}_{\overline{q}}(\eta,\overline{q}\,\beta_1)}{\overline{g}_{\overline{q}}(\overline{q}\,\beta_1,\overline{q}\,\beta_1)} \,\overline{q}\,\beta_1 - \frac{\overline{g}_{\overline{q}}(\eta,\overline{q}'\,\beta_2)}{\overline{g}_{\overline{q}}(\overline{q}'\,\beta_2,\overline{q}'\,\beta_2)} \, \overline{q}'\,\beta_2 \\
	& = \eta - \frac{\langle\eta,\overline{q}\rangle}{\langle\overline{q},\overline{q}\rangle} \, \overline{q} - \frac{\langle\eta,\overline{q}'\rangle}{\langle\overline{q}',\overline{q}'\rangle} \, \overline{q}' .
\end{align*}

The distance between two points $q_1$ and $q_2$ on $\mathcal{Q}_\text{D}$ is defined as
\begin{align*}
	\dist(q_1,q_2) & \eqdef \min_{\sigma\in\Sbb^1} \cos^{-1}\left( \frac{\langle {\overline{q}}_1(\cdot), {\overline{q}}_2(\cdot+\sigma) \rangle}{\norm{\overline{q}_1}_2 \, \norm{\overline{q}_2}_2} \right) \\
	& = \cos^{-1}\left( \frac{\langle {\overline{q}}_1(\cdot), {\overline{q}}_2(\cdot+\sigma_*) \rangle}{\norm{\overline{q}_1}_2 \, \norm{\overline{q}_2}_2} \right) ,
\end{align*}
where $\sigma_*$ denotes the phase shift achieving this minimization. The phase shift $\sigma_*$ corresponds to the phase shift maximizing the circular cross-correlation in \eqref{eq:opt-phase-shift}.

\clearpage
\section{Sensitivity Analysis in the Space of Phase Response Curves} \label{sec:sensitivity}

Sensitivity analysis for oscillators has been widely studied in terms of sensitivity analysis of periodic orbits~\cite{Kramer:1984jk,Rosenwasser:1999tw,Ingalls:2004wr,Wilkins:2009kq}.  
This section develops a sensitivity analysis for phase response curves.  
The sensitivity formula and the developments in this section are closely related to those in~\cite{Vytyaz:2008cw}, which studies the sensitivity analysis of phase response curves, also called perturbation projection vectors, in the context of electronic circuits. The use of sensitivity analysis of phase response curves is novel in the context of biological applications.

This section summarizes the sensitivity analysis for oscillators described by nonlinear time-invariant state-space models with one parameter
\begin{subequations} \label{eq:sys-para}
	\begin{align}
	    \dot{x} & = \fvec(x,u,\lambda), \\
	         y  & = \hvec(x,\lambda),
	\end{align}	
\end{subequations}
where the constant parameter $\lambda$ belongs to some subset $\Lambda\subseteq\Rbb$. 
The scalar nature of the parameter is for convenience but all developments generalize to the multidimensional case.

\subsection{Sensitivity Analysis of a Periodic Orbit}

The (zero-input) steady-state behavior of an oscillator model (that is, its periodic orbit $\perorb$) is characterized by an angular frequency~$\omega(\lambda)$, which measures the speed of a solution along the orbit, and a $2\pi$-periodic steady-state solution $\persol(\cdot;\lambda) = \flow(\cdot/\omega(\lambda),\persol_0(\lambda),\zeroinput,\lambda)$, which describes the locus of this orbit in the state space. 

The sensitivity of the angular frequency at a nominal parameter value $\nlambda$ is the scalar $S^\omega(\nlambda) \in \Rbb$, defined as
\begin{equation*}
	S^\omega(\nlambda) \eqdef \frac{d\omega}{d\lambda}(\nlambda) = \lim_{h \rightarrow 0} \frac{\omega(\nlambda + h) - \omega(\nlambda)}{h} .
\end{equation*}
Likewise, the sensitivity of the $2\pi$-periodic steady-state solution is the $2\pi$-periodic function $S^{\persol}(\cdot;\nlambda):\Sbb^1 \rightarrow \Rbb^{n}$, defined as
\begin{equation*}
	S^{\persol}(\cdot;\nlambda) \eqdef \frac{d\persol}{d\lambda}(\cdot;\nlambda) = \lim_{h \rightarrow 0} \frac{\persol(\cdot;\nlambda + h) - \persol(\cdot;\nlambda)}{h} .
\end{equation*}
From \eqref{eq:bvp_x} and then taking derivatives with respect to $\lambda$,
\begin{subequations} \label{eq:bvp_Zx}
	\begin{align}
			\frac{dS^{\persol}}{d\theta}(\theta;\nlambda) - \frac{1}{\omega} \, A(\theta;\nlambda) \, S^{\persol}(\theta;\nlambda) + \frac{1}{\omega^2} \, v(\theta;\nlambda)\,S^{\omega}(\nlambda) - \frac{1}{\omega} \, E^{\persol}(\theta;\nlambda) & = 0, \\
			S^{\persol}(2\pi;\nlambda) - S^{\persol}(0;\nlambda) & = 0, \\
			\frac{\partial \PLC}{\partial x}(\persol(0;\nlambda);\nlambda) \, S^{\persol}(0;\nlambda) + \frac{\partial \PLC}{\partial \lambda}(\persol(0;\nlambda);\nlambda) & = 0 ,
	\end{align}	
\end{subequations}
where
\begin{align*}
	A(\theta;\nlambda)        & \eqdef \frac{\partial \fvec}{\partial x}(\persol(\theta;\nlambda),0,\nlambda),  \\
	E^{\persol}(\theta;\nlambda)        & \eqdef \frac{\partial \fvec}{\partial \lambda}(\persol(\theta;\nlambda),0,\nlambda),   \\
	v(\theta;\nlambda)& \eqdef \fvec(\persol(\theta;\nlambda),0,\nlambda).
\end{align*}

\begin{rem}
	In the literature~\cite{Ingalls:2004wr,Stelling:2004do,Wilkins:2007is,Bagheri:2007bo,Gunawan:2007jv}, the sensitivity of the period is often preferred to the sensitivity of the angular frequency. It is the real scalar $S^T$
\begin{equation*}
	S^T(\nlambda) \eqdef \frac{dT}{d\lambda}(\nlambda) = \lim_{h \rightarrow 0} \frac{T(\nlambda + h) - T(\nlambda)}{h} .
\end{equation*}	
Both sensitivity measures are equivalent up to a change of sign and a scaling factor, that is, $S^T(\nlambda) / T(\nlambda) = - S^\omega(\nlambda) / \omega(\nlambda)$.
\end{rem}

\subsection{Sensitivity Analysis of a Phase Response Curve}

The input--output behavior of an oscillator model is characterized by its infinitesimal phase response curve $q(\cdot;\lambda)$.

The sensitivity of the infinitesimal phase response curve at a nominal parameter value~$\nlambda$ is the $2\pi$-periodic function $S^{\iPRCu}(\cdot;\nlambda):\Sbb^1\rightarrow\Rbb$, defined as
\begin{equation*}
	S^{\iPRCu}(\cdot;\nlambda) \eqdef \frac{d\iPRCu}{d\lambda}(\cdot;\nlambda) = \lim_{h \rightarrow 0} \frac{\iPRCu(\cdot;\nlambda+h) - \iPRCu(\cdot;\nlambda)}{h}.
\end{equation*}
From \eqref{eq:dirder} and then taking derivatives with respect to $\lambda$,
\begin{align*}
		S^{\iPRCu}(\theta;\nlambda) & = \left\langle S^{\iPRCx}(\theta;\nlambda) , \frac{\partial \fvec}{\partial u}(\persol(\theta;\nlambda),0,\nlambda) \right\rangle \nonumber \\ & \quad + \left\langle \iPRCx(\theta;\nlambda) ,  \frac{\partial^2 \fvec}{\partial x \partial u}(\persol(\theta;\nlambda),0,\nlambda)  \, S^{\persol}(\theta;\nlambda) +  \frac{\partial^2 \fvec}{\partial \lambda \partial u}(\persol(\theta;\nlambda),0,\nlambda) \right\rangle ,
\end{align*}
where the $2\pi$-periodic function $S^{\iPRCx}(\cdot;\nlambda):\Sbb^1\rightarrow\Rbb^{n}$ is the sensitivity of the gradient of the asymptotic phase map evaluated along the periodic orbit $p(\cdot)$, defined as
\begin{equation*}
	S^{\iPRCx}(\cdot;\nlambda) \eqdef \frac{d\iPRCx}{d\lambda}(\cdot;\nlambda) = \lim_{h \rightarrow 0} \frac{\iPRCx(\cdot;\nlambda+h) - \iPRCx(\cdot;\nlambda)}{h}.
\end{equation*}
From \eqref{eq:bvp_q} and then taking derivatives with respect to $\lambda$,
\begin{subequations} \label{eq:bvp_Z_iPRCx}
\begin{align}
		\frac{dS^{\iPRCx}}{d\theta}(\theta;\nlambda) + \frac{1}{\omega} \, \trans{A(\theta;\nlambda)} \, S^{\iPRCx}(\theta;\nlambda) + \frac{1}{\omega} \, \trans{E^\iPRCx(\theta;\nlambda)} \, \iPRCx(\theta;\nlambda) & = 0, \\
		S^{\iPRCx}(2\pi;\nlambda) - S^{\iPRCx}(0;\nlambda) & = 0, \\
		\left\langle S^{\iPRCx}(\theta;\nlambda) , v(\theta;\nlambda) \right\rangle + \left\langle \iPRCx(\theta;\nlambda) , S^v(\theta;\nlambda) \right\rangle - S^\omega(\nlambda) & = 0 ,
\end{align}
\end{subequations}
where 
\begin{align*}
	E^\iPRCx_{ij}(\theta;\nlambda) & \eqdef \sum_{k=1}^{n} \frac{\partial^2 \fvec_i}{\partial x_j \partial x_k}(\persol(\theta;\nlambda),0,\nlambda) \, S^{\persol}_k(\theta;\nlambda) \nonumber \\ & \quad + \frac{\partial^2 \fvec_i}{\partial x_j \partial \lambda}(\persol(\theta;\nlambda),0,\nlambda) - \frac{1}{\omega}\,\frac{\partial \fvec_i}{\partial x_j}(\persol(\theta;\nlambda),0,\nlambda) \, S^\omega(\nlambda) ,\\
	S^v(\theta;\nlambda) & \eqdef \frac{\partial \fvec}{\partial x}(\persol(\theta;\nlambda),0,\nlambda) \, S^{\persol}(\theta;\nlambda) + \frac{\partial \fvec}{\partial \lambda}(\persol(\theta;\nlambda),0,\nlambda) . 
\end{align*}

\subsection{Numerics of Sensitivity Analysis}

Numerical algorithms to solve boundary value problems \eqref{eq:bvp_Zx} and \eqref{eq:bvp_Z_iPRCx} are reviewed in ``Numerical Tools''. Existing algorithms that compute periodic orbits and infinitesimal phase response curves are easily adapted to compute the sensitivity functions of  periodic orbits and   infinitesimal phase response curves, essentially at the same computational cost.

\clearpage
\section{Applications to Biological Systems} \label{sec:illustrations}

This section illustrates the relevance of sensitivity analysis on three system-theoretic case studies arising from biological systems, emphasizing the novel insight provided by the approach described in this article with respect to the existing literature. The first application analyzes the robustness to parameter variations of a circadian oscillator model based on the sensitivity of its phase response curve.  The second application identifies the parameter values of a simple circadian oscillator model in order to fit an experimental-like phase response curve. The third application classifies neural oscillator models based on their phase response curves. All numerical tests were performed with a Matlab numerical code available from the first author's webpage~\cite{Sacre:2008aa}.

\clearpage
\subsection{Robustness Analysis to Parameter Variations: a Case Study in a Quantitative Circadian Oscillator Model}

Testing the robustness of a model against parameter variations is a basic system-theoretic question. In many situations, modeling can be used specifically for the  purpose of identifying the parameters that influence a system property of interest.

In the literature, robustness analysis of circadian rhythms mostly studies the zero-input steady-state behavior (period, amplitude of oscillations, etc.)~\cite{Gonze:2002im,Stelling:2004do,Wilkins:2007is} and empirical phase-based performance measures~\cite{Bagheri:2007bo,Gunawan:2007jv,Taylor:2008co,Hafner:2010us}.

This section defines scalar robustness measures to quantify the sensitivity of the angular frequency (or the period) and the sensitivity of the infinitesimal phase response curve to parameter variations. These robustness measures are applied to a model of the circadian rhythm. A more detailed analysis of this application was presented in \cite{Sacre:2012fk}.

\subsubsection{Scalar Robustness Measure in the Space of Phase Response Curves} 

The angular frequency~$\omega$ is a positive scalar. The sensitivity of~$\omega$ with respect to the parameter~$\lambda$ is thus a real scalar~$S^\omega$, leading to the scalar robustness measure~$R^{\omega} \eqdef \left| S^{\omega} \right|$.
In contrast, the infinitesimal phase response curve (or its equivalence class) $\iPRCu$ belongs to a (nonlinear) space $\mathcal{Q}$. The sensitivity of $\iPRCu$ is thus a vector~$S^{\iPRCu}$ that belongs to the tangent space $\tgspace{\iPRCu}{\mathcal{Q}}$ at $\iPRCu$. A scalar robustness measure~$R^{\iPRCu}$ is defined as
\begin{equation*}
	R^{\iPRCu} \eqdef \left\| S^{\iPRCu} \right\|_{\iPRCu} = \sqrt{g_{\iPRCu}\left(S^{\iPRCu},S^{\iPRCu}\right)},
\end{equation*}
where $\norm{\cdot}_{\iPRCu}$ denotes the norm induced by the Riemannian metric $g_{\iPRCu}\left(\cdot,\cdot\right)$ at ${\iPRCu}$. It is the natural extension of robustness measures to a (nonlinear) space $\mathcal{Q}$. 

When $\mathcal{Q}$ is a quotient space, the element $\iPRCu$ and the tangent vector $S^{\iPRCu}$ are abstract objects. The evaluation of the robustness measure relies on the sensitivity $S^{\overline{\iPRCu}}$ of the signal $\overline{\iPRCu}$ defining the equivalence class in the total space
\begin{equation*}
	R^{\iPRCu} = \left\| P^h_{\overline{\iPRCu}} \, S^{\overline{\iPRCu}} \right\|_{\overline{\iPRCu}} = \sqrt{\overline{g}_{\overline{\iPRCu}}\left(P^h_{\overline{\iPRCu}}\,S^{\overline{\iPRCu}},P^h_{\overline{\iPRCu}}\,S^{\overline{\iPRCu}}\right)},
\end{equation*}
where $P^h_{\overline{\iPRCu}}$ is the projection operator onto the horizontal space $\mathcal{H}_{\bar{q}}$. The projection removes the component of the sensitivity that is tangent to the equivalence class.

When analyzing a model with several parameters ($\lambda \in \Lambda \subseteq \Rbb^l$), all robustness measures~$R^x$ (where $x$ stands for any characteristic of the oscillator) collect the scalar robustness measure corresponding to each parameter in an $l$-dimensional vector.
This vector is often normalized as
\begin{equation*}
	\rho^{x} = \frac{R^x}{\left\|R^x\right\|_{\infty}},
\end{equation*}
where $\norm{\cdot}_{\infty}$ denotes the maximum norm such that components of $\rho^{x}$ belong to the unit interval~$[0,1]$. This measure allows the ranking of model parameters according to their ability to influence the characteristic $x$.

\subsubsection{Quantitative Circadian Oscillator Model} 

\begin{figure}[p]
	\centering
	\includegraphics{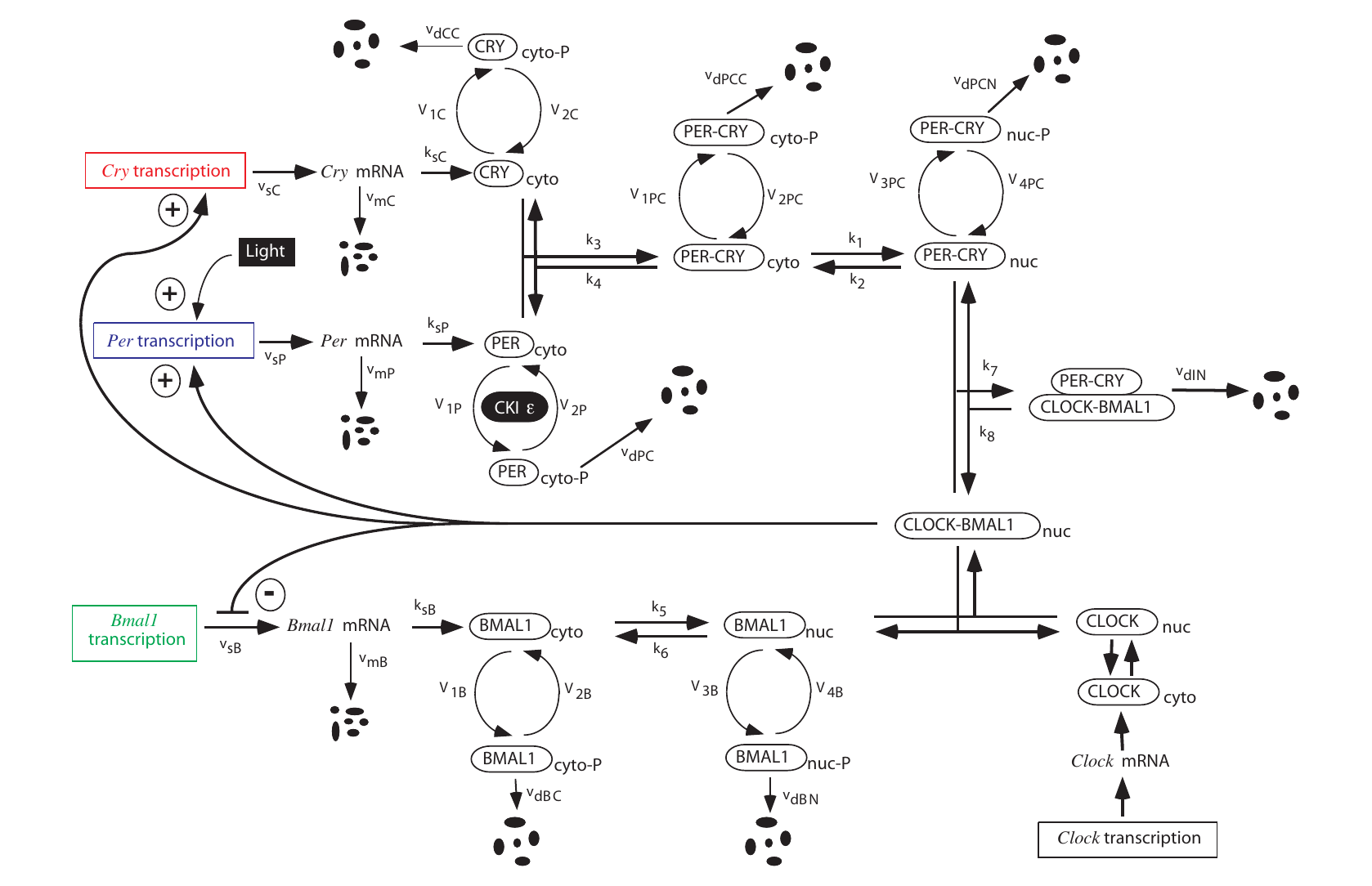}
	\caption{%
		Diagram of the quantitative model for circadian oscillations in mammals involving interlocked negative and positive regulations of \textit{Per}, \textit{Cry}, and \textit{Bmal1} genes by their protein products. (Figure is modified, with permission, from~\cite{Leloup:2003cp}. \copyright{} (2003) National Academy of Sciences, U.S.A.)
	}
	\label{fig:goldbeter-diagram}
\end{figure}

The robustness analysis to parameter variations is illustrated on a quantitative circadian rhythm model for mammals~\cite{Leloup:2003cp}. The model with $16$ state variables and $52$ parameters describes the regulatory interactions between the products of  genes \textit{Per}, \textit{Cry}, and \textit{Bmal1} (see \myfigurename~\ref{fig:goldbeter-diagram}). State-space model equations and nominal parameter values are available in \cite[Supporting Text]{Leloup:2003cp}. The effect of light is incorporated through periodic square-wave variations in the maximal rate of \textit{Per} expression, that is, the value of the parameter $v_{\text{sP}}$ goes from a constant low value during dark phase to a constant high value during light phase. Parameter values remain to be determined experimentally and have been chosen semiarbitrarily within physiological ranges in order to satisfy experimental observations. This model has been extensively studied through unidimensional bifurcation analyses and various numerical simulations of entrainment~\cite{Leloup:2003cp,Leloup:2004co}.

Each parameter of the model describes a single regulatory mechanism, such as transcription and translation control of mRNAs, degradation of mRNAs or proteins, transport reaction, or phosphorylation/dephosphorylation of proteins. 
The analysis of single-parameter sensitivities thus reveals the importance of individual regulatory processes on the function of the oscillator.

In order to enlighten the potential role of circuits rather than single-parameter properties, model parameters were grouped  according to the mRNA loop to which they belonged: \textit{Per}-loop, \textit{Cry}-loop, and \textit{Bmal1}-loop. In addition, parameters associated with interlocked loops were gathered in a last group.

The robustness analysis is developed in the space $\mathcal{Q}_{\text{D}}$ incorporating both scaling and phase-shifting equivalence properties. These equivalence properties are motivated by the uncertainty about the exact magnitude of the light input on the circadian oscillator and by the absence of precise experimental state trajectories, which prevent defining a precise reference position corresponding to the initial phase.

The following section considers sensitivities to relative variations of parameters. 

\subsubsection{Results}

The period and the phase response curve are two characteristics of the circadian oscillator with physiological significance. The sensitivity analysis measures the influence of regulatory processes on tuning the period and shaping the phase response curve.

\begin{figure}[p]
	\centering
	\includegraphics{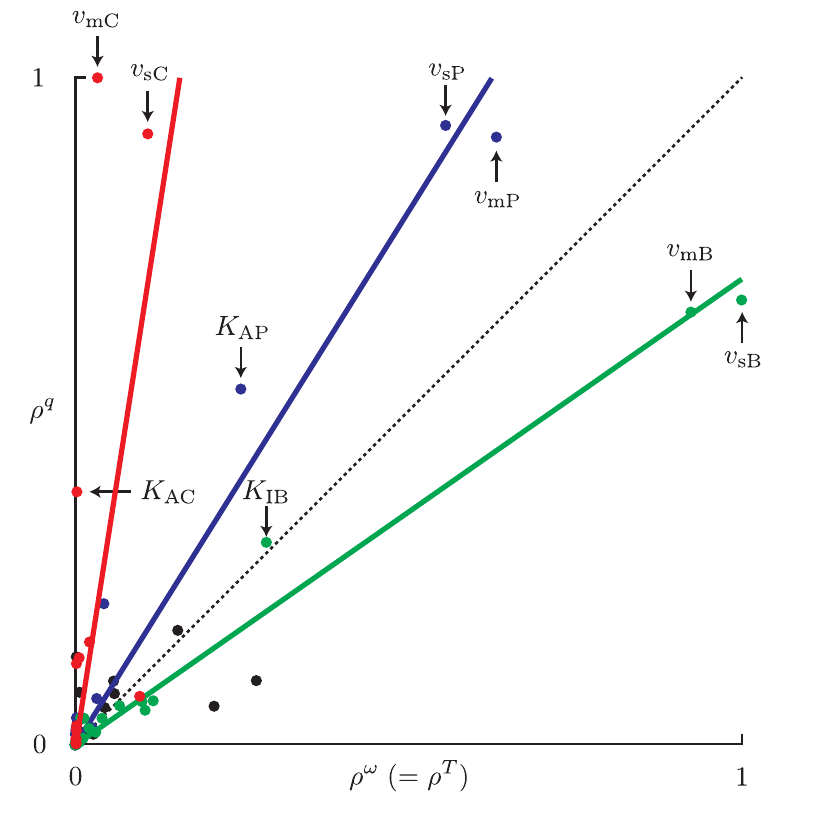}
	\caption{%
		Local robustness analysis to parameter variations in the space of infinitesimal phase response curves.
		Normalized robustness measures $\rho^\omega$ (angular frequency) and $\rho^{\iPRCu}$ (infinitesimal phase response curve) reveal the distinct sensitivity of three distinct genetic circuits (\textit{Cry}, \textit{Per}, and \textit{Bmal1}). 
	Each point is associated with a particular parameter.
	The three lines are regressions over the parameters of the three gene loops.
	The dashed bisector indicates the positions at which the two measures of robustness are identical.
	Only parameters associated with the \textit{Cry}-loop exhibit a low influence on the period and a high influence on the infinitesimal phase response curve.
	The color code corresponds to different subsets of parameters associated with different loops: \textit{Per}-loop in blue, \textit{Cry}-loop in red, and \textit{Bmal1}-loop in green. Parameters associated with interlocked loops are represented in black.
	}
	\label{fig:goldbeter-local}
\end{figure}

A two-dimensional $(\rho^\omega, \rho^{\iPRCu})$ scatter plot in which each point corresponds to a parameter of the model reveals the shape and strength of the relationship between both normalized robustness measures $\rho^\omega$ (angular frequency or, equivalently, period) and $\rho^{\iPRCu}$ (phase response curve). It enables identifying which characteristic is primarily affected by perturbations in individual parameters: parameters below the dashed bisector mostly influence the period; whereas parameters above the dashed bisector mostly influence the phase response curve (see \myfigurename~\ref{fig:goldbeter-local}).

At a coarse level of analysis, the scatter plot reveals that 
the period and the phase response curve exhibit a low sensitivity to most parameters (most points are close to the origin); the period and the phase response curve display a medium or high sensitivity to only few parameters, respectively.

At a finer level of analysis, the scatter plot reveals a qualitative difference in sensitivity to parameters associated with each of the three mRNA loops. The qualitative tendency among parameters associated with the same mRNA loop is represented by a least-square regression line passing through the origin. The following observations are summarized in \mytablename~\ref{tab:goldbeter}.
\begin{itemize}
	\item The \textit{Bmal1}-loop parameters have a strong influence on the period and a medium influence of the phase response curve (regression line below the bisector);
	\item the \textit{Per}-loop parameters have a medium influence on the period  and a high influence on the phase response curve (regression line above the bisector);
	\item the \textit{Cry}-loop parameters have a low influence on the period and a high influence on the phase response curve (regression line above the bisector, close to the vertical axis).
\end{itemize}
In each feedback loop, the three most influential parameters represent the three same biological functions: the maximum rates of mRNA synthesis ($v_\text{sB}$, $v_\text{sP}$, and $v_\text{sC}$), the maximum rate of mRNA degradation ($v_\text{mB}$, $v_\text{mP}$, and $v_\text{mC}$), and the inhibition~(I) or activation~(A) constants for the repression or enhancement of mRNA expression by BMAL1 ($K_\text{IB}$, $K_\text{AP}$, and~$K_\text{AC}$). These three parameters primarily govern the sensitivity associated with each loop.

Two of the three influential parameters of the \textit{Cry}-loop detected by the (local) sensitivity analysis have been identified by numerical simulations as critical for entrainment properties of the model without affecting the period ($K_{\text{AC}}$ in \cite{Leloup:2003cp} and $v_{\text{mC}}$ in \cite{Leloup:2004co}). The local approach supports the importance of these two parameters and identifies the potential importance of a third parameter ($v_{\text{sC}}$).

The conclusions in \cite{Leloup:2003cp} and \cite{Leloup:2004co} rely on extensive simulations of the model under entrainment conditions while varying one parameter at a time. In contrast, the local analysis in this article allows a computationally efficient screening of all parameters. The plot in \myfigurename~\ref{fig:goldbeter-local} was generated in less than a minute with a Matlab code.

\begin{figure}[p!]
	\centering
	\includegraphics{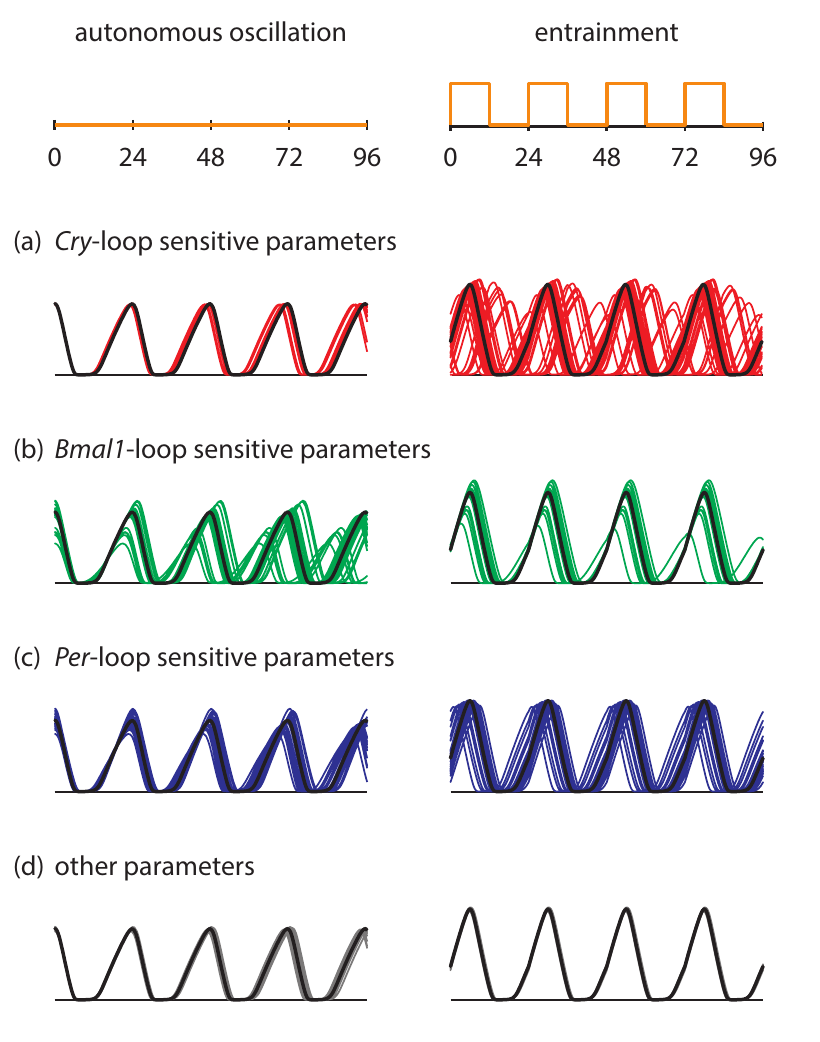}
	\caption{(Caption on next page.)}
	\label{fig:goldbeter-global}
\end{figure}

To evaluate the relevance of the infinitesimal predictions, \myfigurename~\ref{fig:goldbeter-global} displays the time behavior of solutions for different finite parameter changes. The left column illustrates the autonomous solution of the isolated oscillator, and the right column illustrates the steady-state solution entrained by a periodic light input.
Parameter perturbations are randomly taken in a range of $\pm 10\%$ around the nominal parameter value. Each row corresponds to the perturbation of a different group of parameters (the black line corresponds to the nominal system behaviors for nominal parameter values).
\begin{enumerate}[(a)]
	\item 
	Perturbations of the three most influential parameters of the \textit{Cry}-loop ($v_{\text{sC}}$, $v_{\text{mC}}$, and $K_{\text{AC}}$) lead to small  variations (mostly shortening) of the autonomous period and (unstructured) large variations of the phase-locking. This observation is consistent with the low sensitivity of the period and the high sensitivity of the phase response curve. 
	
	\item 
	Perturbations of the three most influential parameters of the \textit{Bmal1}-loop ($v_{\text{sB}}$, $v_{\text{mB}}$, and $K_{\text{IB}}$) lead to medium variations of the autonomous period and phase-locking. The variations of the phase-locking exhibit the same structure as variations of the period, suggesting that the change in period is responsible for the change of phase-locking for these parameters. This observation is consistent with the high sensitivity of the period and the medium sensitivity of the phase response curve. 	
	
	\item 
	Perturbations of the three most influential parameters of the \textit{Per}-loop ($v_{\text{sP}}$, $v_{\text{mP}}$, and $K_{\text{AP}}$) exhibit an intermediate behavior between  situations (a) and (b).
	
	\item 
	Perturbations of parameters of interlocked loops lead to small variations of the autonomous period and the phase-locking, which is consistent with their low sensitivity.
\end{enumerate}
These (nonlocal) observations are thus well predicted by the classification of parameters suggested by the (local) sensitivity analysis (see \myfigurename~\ref{fig:goldbeter-local}).

\clearpage
\subsection{System Identification in the Parameter Space: a Case Study in a Qualitative Circadian Oscillator Model}

System identification builds mathematical models of dynamical systems from observations. In particular, system identification in the parameter space finds a set of parameter values that best match observed data for a given state-space model structure.

Parameter values for circadian rhythm models are often determined by trial-and-error methods due to scant experimental parameter value information.

This section provides a gradient-descent algorithm to identify parameter values that give a phase response curve close to an experimental phase response curve (in a metric described in this article). This algorithm is illustrated on a qualitative circadian oscillator model.
 
\subsubsection{Gradient-Descent Algorithm in the Space of Phase Response Curves} 

A standard technique is to recast the system identification problem as an optimization problem. The parameter estimate is the minimizer of an empirical cost $\tilde{V}(\lambda)$, that is,
\begin{equation*}
	\hat{\lambda} = \argmin_{\lambda\in\Lambda} \tilde{V}(\lambda),
\end{equation*}
where $\tilde{V}(\lambda):\Lambda\rightarrow \Rbb_{\geq0}$ penalizes the discrepancy between observed data and model prediction. Local minimization is usually achieved with a gradient-descent algorithm, requiring the computation of the gradient $\nabla_\lambda \tilde{V}(\lambda)$.

Given an experimental-like phase response curve $\overline{q}_0$ (or its equivalence class $q_0 = [\overline{q}_0]$), a natural cost function $\tilde{V}(\lambda)$ is 
\begin{equation*}
	\tilde{V}(\lambda) \eqdef V(q({\lambda})) = \frac{1}{2} \dist(q(\lambda),q_0)^2,
\end{equation*}
where $\dist(\cdot,\cdot)$ is the distance in the (nonlinear) space $\mathcal{Q}$. 
The gradient (in the parameter space $\Lambda\subseteq\Rbb^l$) of this cost function with respect to the parameter~$\lambda_j$ is
\begin{equation*}
	\nabla_{\lambda_j} \tilde{V}(\lambda) = g_q\left(\grad_q V(q(\lambda)) , S^q_j(\lambda)\right),
\end{equation*}
where $\grad_q V(q(\lambda))$ and $S^q_j(\lambda)$ are elements in the tangent space~$\tgspace{q}{\mathcal{Q}}$.

When $\mathcal{Q}$ is a quotient space, the evaluation of the gradient $\nabla_{\lambda_j} \tilde{V}(\lambda)$ relies on representatives in the total space 
\begin{equation*}
	\nabla_{\lambda_j} \tilde{V}(\lambda) = \overline{g}_{\overline{q}}\left(\grad_{\overline{q}} \overline{V}(\overline{q}(\lambda)) , P^h_{\overline{q}}S^{\overline{q}}_j(\lambda)\right),
\end{equation*}
where $\overline{V}(\overline{q}) = V([\overline{q}])$ for all $\overline{q} \in [\overline{q}]$.

\begin{rem}
	Experimental phase response curves are actually finite discrete sets of measurements. 
	The example problem can be seen as the second step in a procedure in which the first step was fitting a continuous curve to experimental data. 
	This example problem can also be seen as fitting the parameter of a reduced model to reproduce the phase response curve of a detailed, high-dimensional model. In this latter case, the phase response curve of the detailed model serves as the experimental phase response curve.
\end{rem}

\subsubsection{Qualitative Circadian Oscillator Model} 

The system identification is illustrated on a qualitative circadian rhythm model \cite{Goodwin:1965ur}.
The model with 3 state variables and 8 parameters is a cyclic feedback system where metabolites repress the enzymes that are essential for their own synthesis by inhibiting the transcription of the molecule DNA to messenger RNA (see \myfigurename~\ref{fig:goodwin-diagram}). 
It can be described as the cyclic interconnection of three first-order subsystems and a monotonic static nonlinearity
\begin{align*}
	\tau_m\,\dot{m} & = - m + K_m \, \frac{1}{1+[(p+u)/\kappa]^\nu}, \\
	\tau_e\,\dot{e} & = - e + K_e \, m, \\
	\tau_p\,\dot{p} & = - p + K_p \, e. 
\end{align*}
A dimensionless form of this system is equivalent to the constraint  $K_e=K_p=\tau_m=\kappa=1$. For convenience, the remaining static gain is denoted $K_m=K$.

To facilitate the interpretation of the results (but without loss of generality), the parameter space is reduced to two dimensions, by imposing equal time-constants $\tau_e=\tau_p=\tau$ and fixing the Hill coefficient $\nu=20$. 
The parameter space reduces to $(K,\tau)\in\mathbb{R}^2_{>0}$.

The reference phase response curve is chosen in accordance with experimental data and a quantitative circadian rhythm of \textit{Drosophila} \cite{Leloup:1999tx,Leloup:2000kw}. The identification algorithm investigates whether it is possible to match this reference phase response curve with the qualitative Goodwin model. This is a problem for which it is of interest to perform the optimization in the space $\mathcal{Q}_{\text{D}}$ that accounts for scaling and phase-shift invariance. 

\subsubsection{Results}

The Goodwin model exhibits stable oscillations in a region of the reduced parameter space~(see \myfigurename~\ref{fig:simu-exp}a). The border of this region corresponds to a supercritical Andronov-Hopf bifurcation through which the model single equilibrium loses its stability. The contour levels of the cost function---which have been computed in the whole region to make results interpretation easier---reveal two local minima.

Picking initial guess values for model parameters, the gradient-descent algorithm minimizes the cost function following a particular path in the parameter space~(see \myfigurename~\ref{fig:simu-exp}a). The cost function value decreases at each step of the algorithm along this path~(see \myfigurename~\ref{fig:simu-exp}b). The optimal infinitesimal phase response curve (blue or red) is a proper fit for the experimental-like infinitesimal phase response curve (gray), in contrast to the initial infinitesimal phase response curve.

Due to the nonconvexity of the cost function,  paths starting from different initial points may evolve towards different local minima (blue and red paths). In this application, the cost function happens to be (nearly) symmetric with respect to a unitary time-constant~$\tau$, and both local minima correspond to similar infinitesimal phase response curves (up to a scaling factor and a phase shift).

The identification is achieved in the space of infinitesimal phase response curves.
It is of interest to investigate whether the optimal model still compares well to  the quantitative model \cite{Leloup:1999tx}  for non-infinitesimal inputs. \myfigurename~\ref{fig:simu-exp-global} shows that the finite phase response curves of the two models still match. The finite phase response curves were computed through a direct numerical method for the scaling factor and phase shift computed in the optimization procedure.
The shapes of finite phase response curves match. It suggests that finite phase response curves are well captured by the (local) infinitesimal phase response curves.

\clearpage
\subsection{Model Classification in the Parameter Space: a Case Study on a Neural Oscillator Model}
 
Model classification separates models into groups that share common qualitative and/or quantitative characteristics.

Models of neurons are often grouped into two classes based on the bifurcation at the onset of periodic firing~\cite{Hansel:1995vp}. Class-I excitable neuron models exhibit saddle-node-on-invariant-circle bifurcations and can theoretically fire at arbitrarily low finite frequencies. Class-II excitable neuron models exhibit subcritical or supercritical Andronov-Hopf bifurcations and possess a nonzero minimum frequency of firing. Several articles have suggested that class-II neurons display a higher degree of stochastic synchronization than class-I neurons \cite{Galan:2006kx,Galan:2007gi,Galan:2007ha,Marella:2008wj,Abouzeid:2009uh,Hata:2011uq}. All these studies analyze phase models using canonical phase response curves associated with each class (see below) and stress the role played by the shape of the infinitesimal phase response curves for this property. However, the shape of the infinitesimal phase response curve can change quickly once the oscillator model is away from the bifurcation, and thus the qualitative synchronization behavior may also change.

This section compares the usual model classification (class-I versus class-II) to a classification directly based on the distance to canonical infinitesimal phase response curves in the space of phase response curves (class-$q_{\text{I}}$ versus class-$q_{\text{II}}$).

\subsubsection{Model Classification Scheme in the Space of Phase Response Curves} 

A strong relationship between the bifurcation type and the shape of the infinitesimal phase response curve has been demonstrated~\cite{Hansel:1995vp,Ermentrout:1996wr,Brown:2004iy}. Near the bifurcation, the infinitesimal phase response curve of class-I excitable neurons is nonnegative or nonpositive and approximated by
\begin{equation*}
	q_{\text{I}}(\theta) = \left[ 1 - \cos(\theta) \right] ;
\end{equation*}
whereas the infinitesimal phase response curve of class-II excitable neurons has both positive and negative parts and is approximated by
\begin{equation*}
	q_{\text{II}}(\theta) = \sin(\theta + \pi).
\end{equation*}

A model classification based on the distance between the model infinitesimal phase response curve $q$ and canonical infinitesimal phase response curves $q_{\text{I}}$ and $q_{\text{II}}$ is defined as
\begin{equation*}
	q \in
	\begin{cases}
		\text{class-$q_{\text{I}}$}  & \text{if $\dist(q,q_{\text{I}}) < \dist(q,q_{\text{II}})$}, \\
		\text{class-$q_{\text{II}}$} & \text{if $\dist(q,q_{\text{I}}) > \dist(q,q_{\text{II}})$}, \\
	\end{cases} 
\end{equation*}
where $\dist(\cdot,\cdot)$ is the distance in the space $\mathcal{Q}$.

\begin{rem}
	It has been shown that, arbitrarily close to a saddle-node-on-invariant-circle bifurcation, the phase response curve continuously depends on model parameters and its shape can be not only primarily positive or primarily negative but also nearly sinusoidal \cite{Ermentrout:2012cj}.
	However, it remains true that many neural oscillators undergoing a saddle-node-on-invariant-circle bifurcation are such that they exhibit a primarily positive (or primarily negative) phase response curve.
\end{rem}

\subsubsection{Neuron Oscillator Model} 

The model classification is illustrated on a simple two-dimensional reduced model of excitable neurons~\cite{Morris:1981iu}. The model with 2 state variables and 13 parameters is composed of a membrane capacitance in parallel with conductances that depend on both voltage and time (see \myfigurename~\ref{fig:morris-lecar-diagram})  
\begin{align*}
	C\,\dot{V} & = -\overline{g}_{\text{Ca}} \, m_\infty(V) \, (V-V_{\text{Ca}}) - \overline{g}_\text{K} \, w \, (V-V_\text{K}) \nonumber \\ & \quad - \overline{g}_\text{L} \, (V-V_\text{L}) + I_{\text{app}} , \\
	\dot{w} & = \phi \, \frac{w_\infty(V) - w}{\tau_w(V)} ,
\end{align*}
where
\begin{align*}
	m_\infty(V) & = 0.5 \, [1 + \tanh((V-V_1)/V_2)], \\ 
	w_\infty(V) & = 0.5 \, [1 + \tanh((V-V_3)/V_4)], \\
	\intertext{and}
	\tau_w(V) & = 1/\cosh((V-V_3)/(2\,V_4)).
\end{align*}
The applied current $I_{\text{app}}$ is the input.

This model exhibits both classes of excitability for different parameter values \cite{Rinzel:1998ui,Tsumoto:2006wc}. 
For large values of the calcium conductance~$\overline{g}_{\text{Ca}}$, the model exhibits a class-I excitability (saddle-node-on-invariant-circle bifurcation). For smaller values of~$\overline{g}_{\text{Ca}}$, the model exhibits a class-II excitability (Andronov-Hopf bifurcation).

In this context, it is meaningful to classify models based on a distance in the space $\mathcal{Q}_\text{D}$, incorporating both scaling and phase-shifting equivalence properties in order to compare the qualitative shape of infinitesimal phase response curves.

\subsubsection{Results}

The bifurcation-based classification scheme is unidimensional and defines a horizontal separation in the two-dimensional parameter space $(I_{\text{app}},\overline{g}_{\text{Ca}})$ (see \myfigurename~\ref{fig:morris-lecar}a). Indeed, a model is classified based on the bifurcation at the onset of periodic firing while varying the applied current $I_{\text{app}}$.
However, the shape of the infinitesimal phase response curve close to the bifurcation can be different from the canonical shape predicted at the bifurcation boundary (see \myfigurename~\ref{fig:morris-lecar}b).

The classification scheme based on the infinitesimal phase response curve shape provides a different separation in the parameter space (see \myfigurename~\ref{fig:morris-lecar}b). The new classification scheme allows one neuron (for one value of~$\overline{g}_{\text{Ca}}$) to pass from one class to another (crossing of the separation) for different values of applied current~$I_{\text{app}}$. Infinitesimal phase response curves computed for several points close to the bifurcation boundary confirm the classification based on the qualitative shape of infinitesimal phase response curves. In particular, parameter set B belongs to the new class-I.

For class-II oscillators, the correspondence between the bifurcation-based classification and the phase response curve-based classification is limited to a narrow region in the neighborhood of the bifurcation.

To assess the predictive value of the classification, \myfigurename~\ref{fig:morris-lecar-global} displays the time evolution of an uncoupled neuron network in which all neurons are entrained by the same stochastic input (that is, stochastic synchronization). For each neuron (one horizontal line), a point is plotted when the neuron fires (raster plot). Each row (see \myfigurename~\ref{fig:morris-lecar-global}a, from A to C) corresponds to a different set in the parameter space. The synchronization level is quantified by the time evolution of the spike distance in \myfigurename~\ref{fig:morris-lecar-global}b. This distance is equal to 0 for perfect synchronization and to 1 for perfect desynchronization \cite{Kreuz:2013fj}.

The stronger synchronization observed for parameter set C supports the better prediction given by a classification scheme based on the shape of the phase response curve rather than on the bifurcation at the onset of the periodic firing.

\clearpage
\section{Conclusion} \label{sec:conclusion}

The article provided a novel framework to analyze oscillator models in the space of phase response curves and to answer systems questions about oscillator models.
Under some perturbation assumptions, state-space models can be reduced to phase models characterized by their angular frequencies and their phase response curves.

The article proposed to base metrics in the space of dynamical systems on metrics in the space of  phase response curves. 
Quotient Riemannian structures are proposed in order to handle scaling and/or phase invariance properties. The Riemannian framework is used to develop a sensitivity analysis and optimization based analysis or synthesis algorithms.

Three system-theoretic questions arising for biological systems have been considered: robustness analysis to parameter variations, system identification in the parameter space from phase response curve data, and model classification in the parameter space based on distances in the space of phase response curves. While preliminary, these applications suggest that the approach described in this article is numerically efficient and may provide novel insight in several questions of interest for oscillator modeling. 
An inherent limitation of sensitivity analysis is its local nature in the parameter space, in contrast to the global robustness questions encountered in biological applications. The illustrations in ``Robustness Analysis to Parameter Variations: a Case Study in a Quantitative Circadian Oscillator Model'' and elsewhere \cite{Trotta:2012kv} suggest  however that local analyses performed at well chosen operating conditions are good predictors of global trends.

\appendices

\clearpage
\section*{Acknowledgment}
This article presents research results of the Belgian Network DYSCO (Dynamical Systems, Control, and Optimization), funded by the Interuniversity Attraction Poles Programme, initiated by the Belgian State, Science Policy Office. The scientific responsibility rests with its authors. P.~Sacr\'e is a Research Fellow with the Belgian Fund for Scientific Research (F.R.S.-FNRS).

\clearpage


\clearpage

\clearpage

\clearpage

\clearpage

\clearpage
\addtocounter{figure}{-1}
\begin{figure}[p!]
	(Figure on previous page.) 
	\caption{%
		Validation of the local robustness analysis for finite (nonlocal) parameter perturbations.
		Steady-state behaviors for the nominal model and different finite (nonlocal) parameter perturbations are illustrated by time plots of the state variable $M_P$ under constant environmental conditions (autonomous oscillation, left column) and periodic environmental conditions (entrainment, right column).
		 Each row corresponds to the perturbation of a different group of parameters, the black time plot corresponding to the system behavior for nominal parameter values. Perturbations are randomly taken in a range of $\pm 10\%$ around the nominal parameter value (for one parameter at a~time).
		(a)~Perturbations of the three most influential parameters of \textit{Cry}-loop ($v_{\text{sC}}$, $v_{\text{mC}}$, and $K_{\text{AC}}$) lead to small variations of the autonomous period and (unstructured) large variations of the phase-locking.
		(b)~Perturbations of the three most influential parameters of \textit{Bmal1}-loop ($v_{\text{sB}}$, $v_{\text{mB}}$, and $K_{\text{IB}}$) lead to larger variations of the autonomous period and medium variations of the phase-locking.
		(c)~Perturbations of the three most influential parameters of \textit{Per}-loop ($v_{\text{sP}}$, $v_{\text{mP}}$, and $K_{\text{AP}}$) exhibit an intermediate behavior between the situations (a) and (b).
		(d)~Perturbations of parameters of interlocked loops lead to small variations of the autonomous period and the phase-locking.
	}
\end{figure}

\clearpage
\begin{figure}[p]
	\centering
	\includegraphics{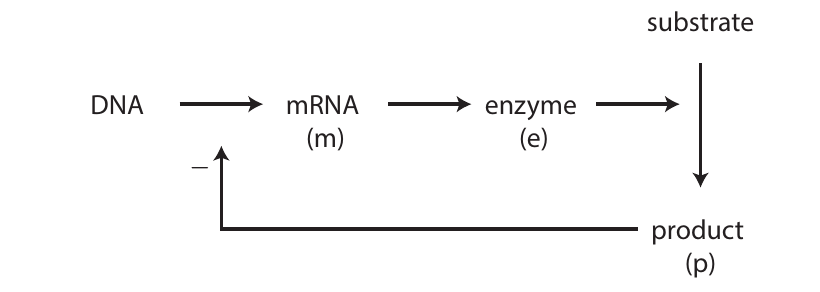}
	\caption{%
		Diagram of the qualitative model for circadian oscillations~\cite{Goodwin:1965ur}. The qualitative model represents the effect of products (p) that repress the enzymes (e) which are essential for their own synthesis by inhibiting the transcription of the molecule DNA to messenger RNA (m).
	}
	\label{fig:goodwin-diagram}
\end{figure}

\clearpage
\begin{figure}[p]
	\centering
	\includegraphics{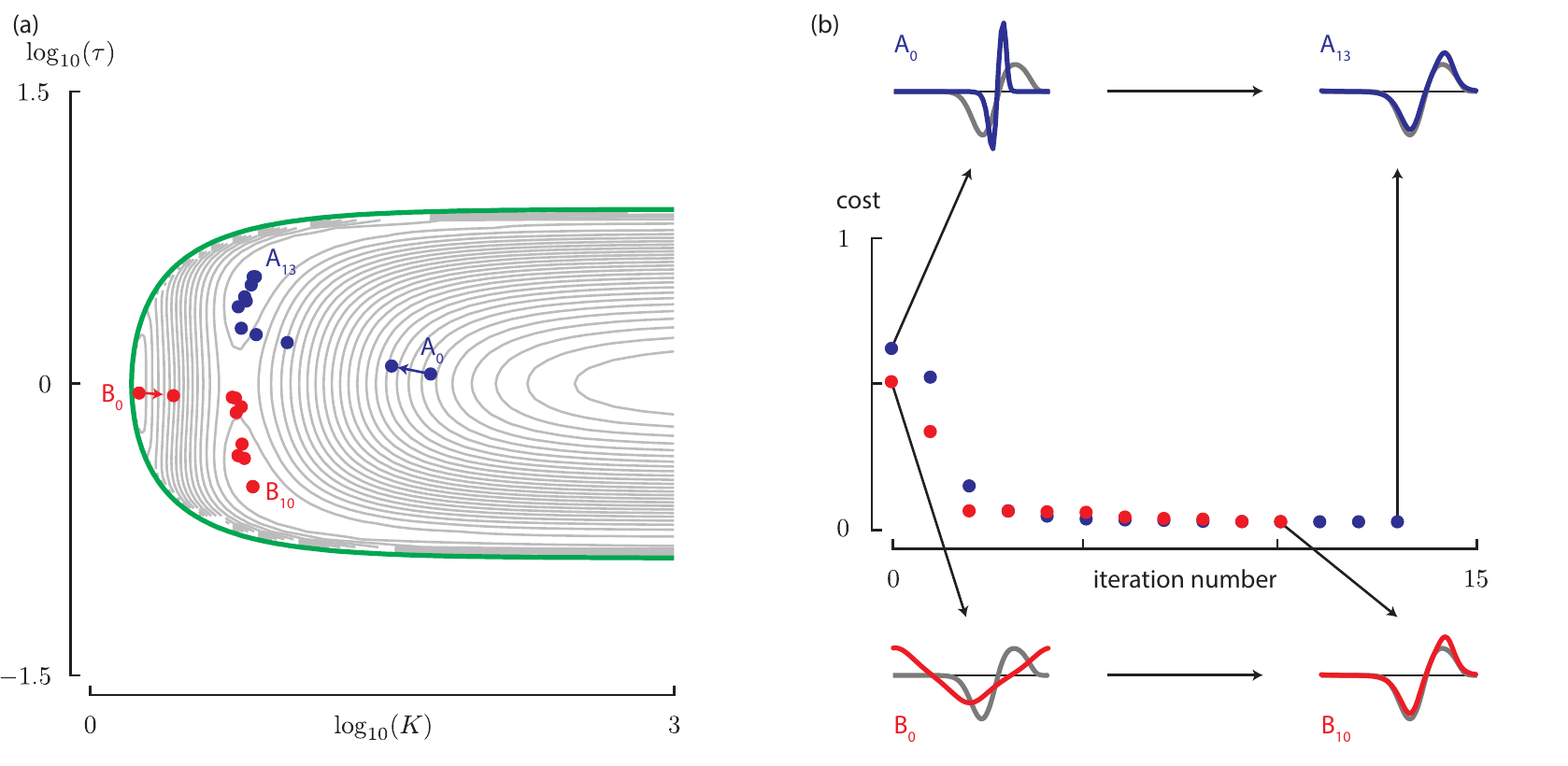}
	\caption{%
		System identification in the parameter space from phase response curve data.
		(a)~The cost function (gray contours) between an experimental infinitesimal phase response curve and the infinitesimal phase response curves exhibits a nonconvex behavior in the reduced parameter space. The gradient-descent algorithm follows the path indicated by dots (two random trials from different initial parameter values are shown in blue and red, respectively). 
		(b)~The cost along the path followed by the gradient-descent algorithm decreases with the iteration number. The shape of the optimal infinitesimal phase response curve (blue or red) is closer to the reference infinitesimal phase response curve (gray) than the initial infinitesimal phase response curve (blue or red).%
		}
	\label{fig:simu-exp}
\end{figure}

\clearpage
\begin{figure}[p]
	\centering
	\includegraphics{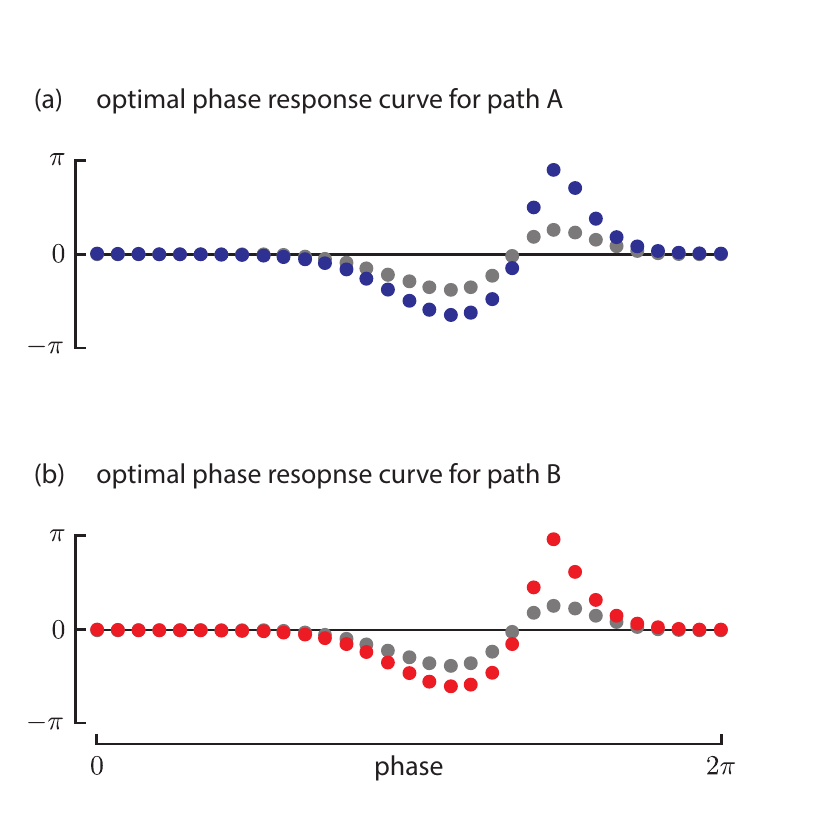}
	\caption{%
		Validation of the system identification for finite phase response curves.
		The finite phase response curves computed at optimally identified parameters (blue or red) in the parameter space match well with the finite phase response curve of the quantitative circadian rhythm model (gray). The magnitude of the input and the reference point have been chosen based on the results of the optimization procedure in the space of infinitesimal phase response curves.
		Subfigures (a) and (b) correspond to the result of the two random trials, respectively.
	}
	\label{fig:simu-exp-global}
\end{figure}

\clearpage
\begin{figure}[p]
	\centering
	\includegraphics{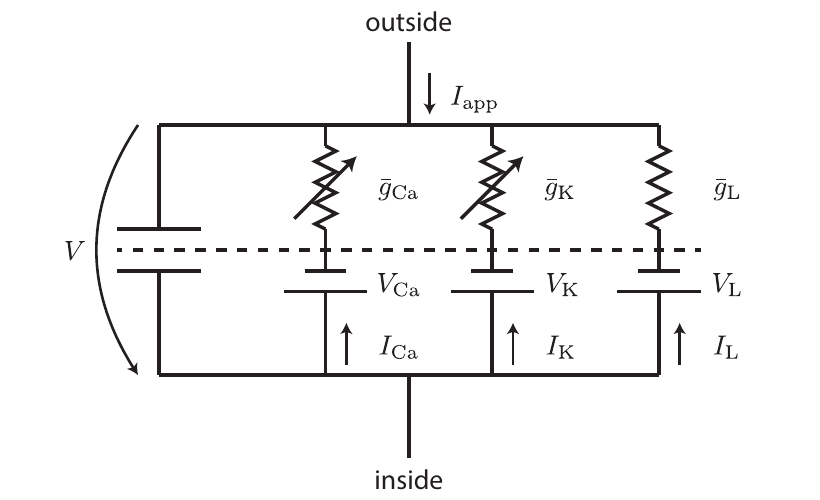}
	\caption{%
		Equivalent circuit diagram of the model for excitable neurons~\cite{Morris:1981iu}. The model is composed of one compartment containing the conductances shown, in parallel with a membrane capacitance.
	}
	\label{fig:morris-lecar-diagram}
\end{figure}

\clearpage
\begin{figure}[p]
	\centering
	\includegraphics{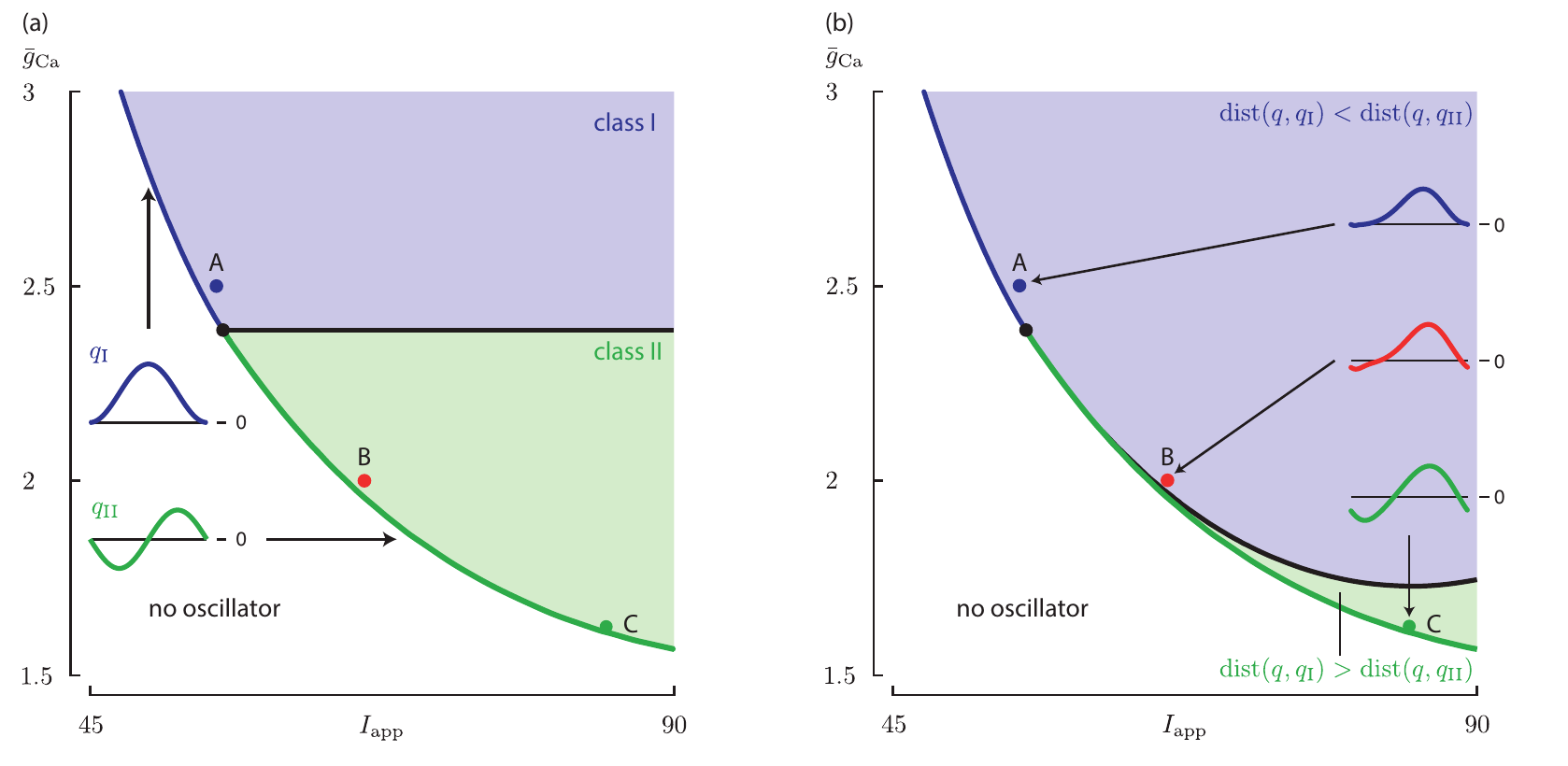}
	\caption{%
		Model classification in the parameter space based on a distance in the space of infinitesimal phase response curves.
		(a)~Standard classification relies on the bifurcation at the onset of the periodic orbit while varying the applied current $I_{\text{app}}$ (class-I in blue and class-II in green). This unidimensional classification defines a horizontal separation in the parameter space. Ideal phase response curves at the bifurcation are shown.
		(b)~The classification relies on the distance to nearest ideal phase response curves (class-I in blue and class-II in green). This classification in the two-dimensional parameter space determines different subsets.
		Parameter set A (respectively C) belongs to class-I (respectively class-II) and its phase response curve is closest to the canonical phase response curve $\iPRCu_{\text{I}}$ (respectively canonical phase response curve $\iPRCu_{\text{II}}$). However, parameter set B (in red) belongs to class-II and its phase response curve is closest to the canonical phase response curve $\iPRCu_{\text{I}}$.
		(Parameter values: $C=20 \; \mu\text{F}/\text{cm}^2$, $\overline{g}_\text{K}=8\;\text{mS}/\text{cm}^2$, $\overline{g}_\text{L}=2\;\text{mS}/\text{cm}^2$, $V_\text{Ca} = 120 \;\text{mV}$, $V_\text{K}=-80 \;\text{mV}$, $V_\text{L}=-60 \;\text{mV}$, $V_1=-1.2 \;\text{mV}$, $V_2=18 \;\text{mV}$, $V_3 = 12\;\text{mV}$, $V_4 = 17.4 \;\text{mV}$, $\phi = 1/15 \;\text{s}^{-1}$.)%
		}
	\label{fig:morris-lecar}
\end{figure}

\clearpage
\begin{figure}[p]
	\centering
	\includegraphics{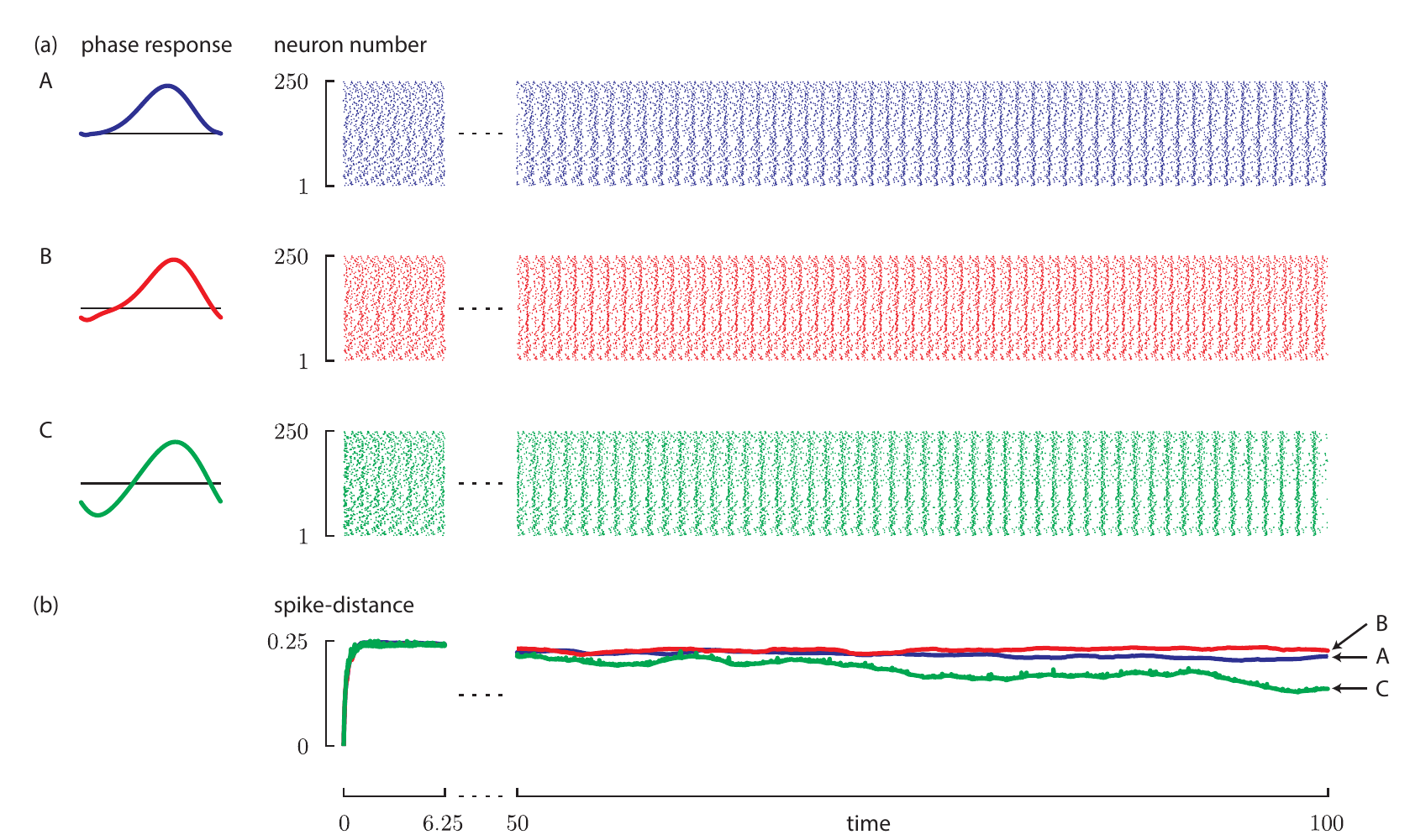}
	\caption{
		Validation of the model classification for  stochastic synchronization.
		(a)~Stochastic synchronization for uncoupled network of state-space models are illustrated by time plots of firing times (one line corresponds to one neuron). Each row (from A to C) corresponds to a point in the parameter space (see \myfigurename~\ref{fig:morris-lecar}). 
		(b)~The spike-distance quantifies the synchronization level of the network (it is equal to 0 for perfect synchronization and to 1 for perfect desynchronization). 
		Parameter sets A and B exhibit a lower stochastic synchronization (higher values of the spike-distance) than parameter set C, consistent with the fact that the phase response curve of parameter set B is shapewise closer to the phase response curve of parameter set A than to the one of parameter set C.
		}
	\label{fig:morris-lecar-global}
\end{figure}

\clearpage
\begin{table}[p]
	\centering
	
	\begin{tabular}{|c||c|c|}
		\hline
		& $q(\cdot) \not\sim q(\cdot)\,\alpha$ & $q(\cdot) \sim q(\cdot)\,\alpha$ \\
		\hline \hline
		$q(\cdot) \not\sim q(\cdot+\sigma)$ & $\mathcal{Q}_\text{A}\eqdef\Hilbert^{1}$ &  $\mathcal{Q}_\text{B}\eqdef\Hilbert^{1}/\Rbb_{>0}$ \\
		$q(\cdot) \sim q(\cdot+\sigma)$ & $\mathcal{Q}_\text{C}\eqdef\Hilbert^{1}/\Shift(\Sbb^1)$ & $\mathcal{Q}_\text{D}\eqdef\Hilbert^{1}/(\Shift(\Sbb^1)\times\Rbb_{>0})$ \\
		\hline
	\end{tabular}
	
	\caption{Combining equivalence properties defines different quotient spaces for phase response curves.}
	\label{tab:metric}
\end{table}

\clearpage
\begin{table}
	\centering
	\begin{tabular}{|c||cc|}
		\hline
		& sensitivity of   & sensitivity of \\
				& the period   & the phase response curve \\

		\hline
		\hline
		\textit{Per} loop & medium  & high \\
		\textit{Cry} loop & low  & high \\
		\textit{Bmal1} loop & high  & medium \\
		\hline
	\end{tabular}
	\caption{The robustness analysis reveals a qualitative difference in sensitivity to parameters associated with each of the three mRNA loops.}
	\label{tab:goldbeter}
\end{table}

\clearpage
\begin{table}[p]
	\centering
	\begin{tabular}{|c||c|c|}
		\hline 
		& Forward Multiple Shooting & Trapezoidal Scheme  \\
		\hline 
		$r_i$ & $\phi\left(\frac{h_i}{\omega},\persol_i,\zeroinput,\lambda\right) - \persol_{i+1}$ & $\persol_{i+1} - \persol_i - \frac{1}{2}\frac{h_i}{\omega}\left[f\left(\persol_{i},0,\lambda\right) + f\left(\persol_{i+1},0,\lambda\right)\right]$ \\[12pt]

		$G_i$ & $\frac{\partial \phi}{\partial x_0}\left(\frac{h_i}{\omega},\persol_i,\zeroinput,\lambda\right)$ & $- I_n - \frac{1}{2}\frac{h_i}{\omega}  \frac{\partial f}{\partial x}\left(\persol_{i},0,\lambda\right)$   \\[3pt]
		$H_i$ & $I_n$ & $- I_n + \frac{1}{2} \frac{h_i}{\omega} \frac{\partial f}{\partial x}\left(\persol_{i+1},0,\lambda\right)$   \\[3pt]
		$b_i^{\persol}$ & $-\frac{h_i}{\omega^2}\frac{\partial \phi}{\partial t}\left(\frac{h_i}{\omega},\persol_i,\zeroinput,\lambda\right)$ & $-\frac{1}{2}\frac{h_i}{\omega^2}\left[f\left(\persol_i,0,\lambda\right) + f\left(\persol_{i+1},0,\lambda\right)\right]$   \\[12pt]

		$\tilde{G}_i$ & $I_n$ & $- I_n + \frac{1}{2}  \frac{h_i}{\omega} \trans{\frac{\partial f}{\partial x}\left(\persol_i,0,\lambda\right)}$   \\[3pt]
		$\tilde{H}_i$ & $\trans{\frac{\partial \phi}{\partial x_0}\left(\frac{h_i}{\omega},\persol_i,\zeroinput,\lambda\right)}$ & $- I_n - \frac{1}{2} \frac{h_i}{\omega} \trans{\frac{\partial f}{\partial x}\left(\persol_{i+1},0,\lambda\right)}$   \\[12pt]
		$E_i^{\persol}$ & $\frac{\partial \phi}{\partial \lambda}\left(\frac{h_i}{\omega},\persol_i,\zeroinput,\lambda\right)$ & $\frac{1}{2}\frac{h_i}{\omega}\left[ E^{\persol}\left(\theta_i;\lambda\right) + E^{\persol}\left(\theta_{i+1};\lambda\right) \right]$ \\
		$E_i^{\iPRCx}$ & $\trans{\left[\frac{d}{d\lambda}\frac{\partial \phi}{\partial x_0}\left(\frac{h_i}{\omega},\persol_i,\zeroinput,\lambda\right)\right]} \iPRCx_{i+1}$ & $-\frac{1}{2}\frac{h_i}{\omega}\left[ \trans{E^{\iPRCx}\left(\theta_i;\lambda\right)} \iPRCx_{i} + \trans{E^{\iPRCx}\left(\theta_{i+1};\lambda\right)} \iPRCx_{i+1} \right]$ \\
		$P$ & $\frac{1}{N+1}I_{\left(N+1\right)n}$ & $\frac{1}{2\pi}\diag\left( \frac{h_0}{2},\frac{h_0+h_1}{2},\ldots,\frac{h_{N-1}+h_{N}}{2},\frac{h_N}{2}\right) \otimes I_n$ \\
		\hline 
	\end{tabular}
	\caption{Residual maps $r_i(\persol_\mesh,\omega)$; linear block entries $G_i$, $H_i$, and $b_i^{\persol}$; adjoint linear block entries $\tilde{G}_i$ and $\tilde{H}_i$; sensitivity block entries $E_i^{\persol}$, $E_i^{\iPRCx}$; and $P$ for two one-step numerical algorithms ($i=0,1,\ldots,N-1$).}
	\label{tab:num_tools}
\end{table}

\setcounter{equation}{0}
\setcounter{figure}{0}
\setcounter{table}{0} 

\renewcommand{\theequation}{S\arabic{equation}}
\renewcommand{\thefigure}{S\arabic{figure}}
\renewcommand{\thetable}{S\arabic{table}} 

\newcounter{sidebar}
\setcounter{sidebar}{1}

\clearpage

\renewcommand{\nomname}{Sidebar \thesidebar: List of Symbols} 
\printnomenclature

\clearpage

\addtocounter{sidebar}{1}
\section*{Sidebar \thesidebar: A Brief History of Phase Response Curves}

Phase response curves were used for the first time in 1960 by a biological experimentalist in order to represent the results of phase-resetting experiments on the rhythm of the daily locomotor activity in flying squirrels~\cite{DeCoursey:1960cw}. The author was investigating the effect of short light pulses on the daily onsets of running activity in the wheel in constant darkness. The response to these stimulations varied according to the time of the day---the squirrel's subjective day---at which the light pulse was administered. To represent her data, the author plotted the observed time shift (advance or delay) as a function of time of perturbation.

\subsection*{Phase Response Curves in Biology}

Phase response curves are widely used to study biological rhythms (see the pioneering book~\cite{Winfree:1980ue}). The two main applications are circadian rhythms and neural (or cardiac) excitable cells.

In circadian rhythms, the phase response curve is used to study the effect of light (and sometimes the effect of drugs, such as melatonin) on the rhythm. In particular, the mechanism of entrainment to light is of critical importance in circadian rhythm studies.
Numerous experimental phase response curves for circadian rhythms have been compiled in an atlas \cite{Johnson:1990vx}. Most of these phase response curves have a typical shape including a dead-zone, which is an interval of zero sensitivity during the subjective day of the studied organism.

In neural (or cardiac) excitable cells, the phase response curve is used to study ensemble behavior in a network: particularly, synchronization in coupled neurons and entrainment in uncoupled neurons subject to correlated inputs (also known as stochastic synchronization). The 2012 book \cite{Schultheiss:2012bz} compiles several applications of phase response curves in neuroscience.

\subsection*{Phase Response Curves in Engineering}

Phase response curves are not often used in engineering applications.
An exception is in electronic circuits, where the concept of perturbation projection vector was developed to study phase noise in oscillators \cite{Kaertner:1989fk,Kaertner:1990kj,Demir:2000kv,Demir:2002kn,Demir:2006kf}. 
Mathematically, the perturbation projection vector is identical to the infinitesimal phase response curve. It is used as a reduction tool to study oscillators and to design electronic circuits \cite{Vytyaz:2009fa}.

\clearpage

\addtocounter{sidebar}{1}
\section*{Sidebar \thesidebar: Phase Maps}

\begin{figure}[p]
	\centering
	\includegraphics{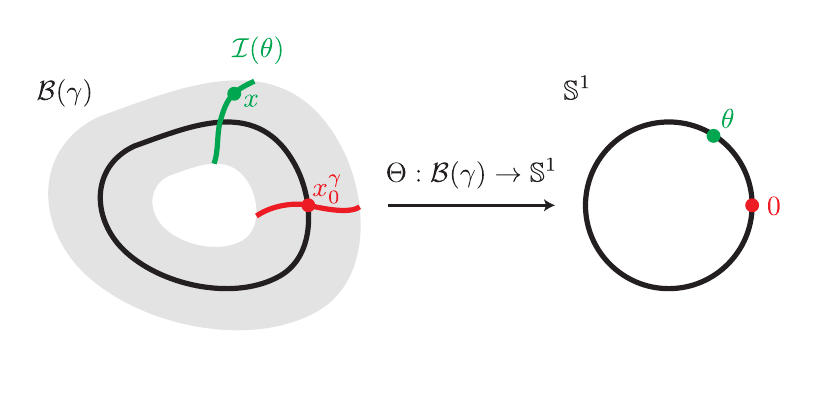}
	\caption{%
		Asymptotic phase map and isochrons.
		The asymptotic phase map $\Theta:\mathcal{B}(\perorb)\rightarrow\Sbb^1$ associates with each point $x$ in the basin of attraction $\mathcal{B}(\perorb)$ a scalar phase $\Theta(x)=\theta$ on the unit circle $\Sbb^1$ such that $\lim_{t\rightarrow+\infty} \norm{\flow(t,x,\zeroinput) - \flow(t + \theta/\omega,\persol_0,\zeroinput)}_2 = 0$. The image of $\persol_0$ through the phase map $\Theta$ is equal to 0. 
		The set of points associated with the same phase $\theta$ (that is, a level set of the phase map) is called an isochron and is denoted by $\mathcal{I}(\theta)$.
		}
	\label{fig:notations}
\end{figure}

Phase maps, as well as the associated notion of isochrons, are key ingredients for studying oscillator models. The brief exposition of phase maps given below follows the terminology and definitions of \cite{Winfree:1980ue} and \cite{Glass:1988ub}. The notation is illustrated in \myfigurename~\ref{fig:notations}.

Consider an oscillator described by \eqref{eq:nlsys}. 

The basin of attraction of $\perorb$ (the oscillator stable set) is the maximal open set from which the periodic orbit $\perorb$ attracts, that is,
\begin{equation*}
	\mathcal{B}(\perorb)\eqdef\{x_0\in\mathcal{X}:\lim_{t\rightarrow+\infty} \dist(\flow(t,x_0,\zeroinput),\perorb) = 0\}
\end{equation*}
where $\dist(x,\perorb) \eqdef \inf_{y\in\perorb} \norm{x - y}_2$ is the distance from the point $x\in\mathcal{X}$ to the set $\perorb\subseteq\mathcal{X}$ based on the Euclidean norm $\norm{ \cdot }_2$ in $\mathbb{R}^n$. 

Since the periodic orbit~$\perorb$ is a one-dimensional manifold in~$\Rbb^n$, it is homeomorphic to the unit circle~$\mathbb{S}^1$. It is thus naturally parameterized in terms of a single scalar phase.
The smooth bijective phase map $\varTheta:\perorb\rightarrow\Sbb^1$ associates with each point $x$ on the periodic orbit $\perorb$ its phase $\varTheta(x)\reveqdef\vartheta$ on the unit circle $\Sbb^1$, such that
\begin{equation*}
	x - \flow(\vartheta/\omega,\persol_0,\zeroinput) = 0.
\end{equation*}
This mapping is constructed such that the image of the reference point $\persol_0$ is equal to~$0$ (that is, $\varTheta(\persol_0)=0$) and the progression along the periodic orbit (in absence of perturbations) produces a constant increase in $\vartheta$.
The phase variable $\vartheta:\Rbb_{\geq 0} \rightarrow \Sbb^1$ is defined along each zero-input trajectory $\flow(\cdot,x_0,\zeroinput)$ starting from a point $x_0$ on the periodic orbit $\perorb$, as $\vartheta(t) \eqdef \varTheta(\flow(t,x_0,\zeroinput))$ for all times $t \geq 0$. The phase dynamics are thus given by $\dot{\vartheta} = \omega$.

For a hyperbolic stable periodic orbit, the notion of phase can be extended to any point~$x$ in the basin~$\mathcal{B}(\perorb)$ by defining the concept of asymptotic phase. 
The asymptotic phase map~$\Theta : \mathcal{B}(\perorb) \rightarrow \mathbb{S}^1$ associates with each point~$x$ in the basin~$\mathcal{B}(\perorb)$ its asymptotic phase~$\Theta(x) \reveqdef \theta$ on the unit circle $\Sbb^1$, such that 
\begin{equation*}
	\lim_{t\rightarrow+\infty} \left\|\flow(t,x,\zeroinput) - \flow(t,\flow(\theta/\omega,\persol_0,\zeroinput),\zeroinput)\right\|_2 = 0.
\end{equation*}
Again, this mapping is constructed such that the image of~$\persol_0$ is equal to~$0$ and such that the progression along any orbit in~$\mathcal{B}(\perorb)$ (in absence of perturbations) produces a constant increase in~$\theta$.
The asymptotic phase variable $\theta:\Rbb_{\geq 0} \rightarrow \Sbb^1$ is defined along each zero-input trajectory $\flow(\cdot,x_0,\zeroinput)$ starting from a point $x_0$ in the basin of attraction of~$\perorb$ as $\theta(t) \eqdef \Theta(\flow(t,x_0,\zeroinput))$ for all times $t \geq 0$. The asymptotic phase dynamics are thus given by $\dot{\theta} = \omega$.

The notion of the asymptotic phase variable can be extended to any nonzero-input trajectory $\flow(\cdot,x_0,u(\cdot))$ in the basin of attraction of $\perorb$. In this case, the asymptotic phase variable is defined as $ \theta(t) \eqdef \Theta(\flow(t,x_0,u(\cdot)))$ for all times $t \ge 0$. 
Thus, the phase variable $\theta(t_*)$, at an instant $t_* \ge 0$, evaluates the asymptotic phase of the point $\flow(t_*,x_0 ,u(\cdot))$. 
The asymptotic phase dynamics in the case of a nonzero input are often hard to derive. 

Level sets of the asymptotic phase map $\Theta$, that is, sets of all points in the basin of $\perorb$ with the same asymptotic phase, are termed isochrons. Formally, the isochron $\mathcal{I}(\theta)$ associated with the asymptotic phase $\theta$ is the set $\mathcal{I}(\theta) \eqdef \left\{ x\in \mathcal{B}(\perorb) : \Theta(x) = \theta \right\}$. Considering hyperbolic periodic orbits, isochrons are codimension-1 submanifolds (diffeomorphic to $\Rbb^{n-1}$) crossing the periodic orbit transversally and foliating the entire basin of attraction~\cite{Guckenheimer:1975tt}.

In general, the (asymptotic) phase maps and their isochrons are complex. This often makes analytical computation impossible and even numerical computation intractable (or at least expensive, particularly for high-dimensional oscillator models). Most numerical techniques rely on backward integration \cite{Guillamon:2009gy,Osinga:2010eu,Sherwood:2010ck}. An elegant forward integration method was developed in \cite{Mauroy:2012vl} and extended to stable fixed points in \cite{Mauroy:2013kx}.

\clearpage

\addtocounter{sidebar}{1}
\section*{Sidebar \thesidebar: From Infinitesimal to Finite Phase Response Curves}

The concept of infinitesimal and finite phase response curves are closely related under the assumption of weak input. Below, the brief exposition highlights the relationship between these two concepts.

By definition, the finite phase response curve $\PRC(\theta;u(\cdot))$ measures the asymptotic difference between the images through the asymptotic phase map $\Theta$ of the perturbed trajectory $\flow(t,\persol(\theta),u(\cdot)))$ and the unperturbed trajectory $\flow(t,\persol(\theta),\zeroinput))$, that is,
\begin{equation} \label{eq:prc_def_supp}
	\PRC(\theta;u(\cdot)) = \lim_{t\rightarrow\infty} [ \Theta(\flow(t,\persol(\theta),u(\cdot))) - \Theta(\flow(t,\persol(\theta),\zeroinput))] \pwrap{[-\pi,\pi)}.
\end{equation}
Linearizing \eqref{eq:prc_def_supp} around the unperturbed trajectory $(\flow^*(t),u^*(t)) \eqdef (\flow(t,\persol(\theta),\zeroinput),\zeroinput)$ and defining the perturbations $(\delta \flow(t),\delta u(t)) \eqdef (\flow(t,\persol(\theta),u(\cdot)) - \flow^*(t),u(t) - u^*(t))$ lead to
\begin{align*}
	\PRC(\theta;u(\cdot)) & = \lim_{t\rightarrow\infty} [ \Theta(\flow^*(t) + \delta\flow(t)) - \Theta(\flow^*(t))] \\
	& = \lim_{t\rightarrow\infty} [ \Theta(\flow^*(t)) + \trans{\nabla_x\Theta(\flow^*(t))}\,\delta \flow(t)  - \Theta(\flow^*(t)) + \mathcal{O}(\norm{\delta \flow(t)}_2^2)] \\
	& = \lim_{t\rightarrow\infty} \trans{\nabla_x\Theta(\flow^*(t))}\,\delta \flow(t) + \mathcal{O}(\norm{\delta \flow(t)}_2^2),
\end{align*}
where the perturbation $\delta\flow(t)$ is the solution of the linearized system
\begin{equation*}
	\delta\dot{\flow}(t) = \underbrace{\frac{\partial \fvec}{\partial x}(\flow^*(t),u^*(t))}_{\reveqdef A_{\flow}(t) = A(\omega\,t+\theta)} \, \delta \flow(t) + \underbrace{\frac{\partial \fvec}{\partial u}(\flow^*(t),u^*(t))}_{\reveqdef b_{\flow}(t) = b(\omega\,t+\theta)} \, \delta u(t) + \mathcal{O}(\norm{\delta \flow}_2^2,\abs{\delta u}^2,\norm{\delta \flow}_2\,\abs{\delta u}).
\end{equation*}
The solution of the linearized equation is 
\begin{equation*}
	\delta \flow(t) = \Phi(t,0)\, \delta \flow(0) + \int_{0}^{t} \Phi(t,s) \, b_{\flow}(s) \, \delta u(s) \,  ds,
\end{equation*}
where the fundamental solution $\Phi(\tau,\sigma)$ associated with $A_\flow(t)$ is the solution of the following matrix equation
\begin{equation*}
	\frac{\partial \Phi}{\partial \tau}(\tau,\sigma) = A_\flow(\tau) \, \Phi(\tau,\sigma), \quad \Phi(\sigma,\sigma) = I_n.
\end{equation*}
The gradient of the asymptotic phase map evaluated along the unperturbed trajectory is given by $\nabla_x\Theta(\flow^*(t)) = \iPRCx(\omega\,t+\theta)$ and is the solution of the adjoint linearized equation~\eqref{eq:bvp_q}. 
Exploiting the properties of the fundamental solution leads to $\trans{\iPRCx(\omega\,t+\theta)}\,\Phi(t,s) = \trans{\iPRCx(\omega\,s+\theta)}$.
Because $\delta\flow(0) = 0$ and $\delta u(t) = u(t)$, we have thus 
\begin{align*}
	\PRC(\theta;u(\cdot)) & \approx \lim_{t\rightarrow\infty} \trans{\iPRCx(\omega\,t+\theta)}\,\left[ \Phi(t,0)\, \delta \flow(0) + \int_{0}^{t} \Phi(t,s) \, b(\omega\,s+\theta) \, \delta u(s) \,  ds \right] \\
		& = \lim_{t\rightarrow\infty}  \int_{0}^{t} \trans{\iPRCx(\omega\,s+\theta)} \, b(\omega\,s+\theta) \, u(s) \,  ds .
\end{align*}
Finally, the finite phase response curve is thus approximated by the ``convolution'' between the infinitesimal phase response curve and the phase-resetting input $u(t)$, that is,
\begin{equation*}
\PRC(\theta;u(\cdot)) \approx \lim_{t\rightarrow\infty}  \int_{0}^{t} \iPRCu(\omega\,s+\theta)\, u(s) \,  ds .
\end{equation*}

\clearpage

\addtocounter{sidebar}{1}
\section*{Sidebar \thesidebar: Basic Concepts of Differential Geometry on Manifolds}

This brief exposition recalls basic concepts of differential geometry on manifolds. It follows the terminology and definitions of \cite{Absil:2008uy}.

A manifold $\mathcal{M}$ is endowed with a Riemannian metric $g_x(\xi_x,\zeta_x)$, which is an inner product of two elements $\xi_x$ and $\zeta_x$ of the tangent space $\tgspace{x}{\mathcal{M}}$ at $x$. The metric induces a norm on $\tgspace{x}{\mathcal{M}}$ at $x$
\begin{equation*}
	\|\xi_x\|_x \eqdef \sqrt{g_x(\xi_x,\xi_x)}.
\end{equation*}

The length of a curve $\gamma:(a,b)\subset\Rbb\rightarrow\mathcal{M}$ is defined as
\begin{equation*}
	L(\gamma) \eqdef \int_{a}^{b} \|\dot{\gamma}(t)\|_{\gamma(t)} dt.
\end{equation*}
The geodesic distance between two points $x$ and $y$ on $\mathcal{M}$ is defined as
\begin{equation*}
	\dist(x,y) = \min_{\Gamma} L(\gamma),
\end{equation*}
where $\Gamma$ is the set of all curves in $\mathcal{M}$ joining points $x$ and $y$
\begin{equation*}
	\Gamma = \{ \gamma:[0,1]\rightarrow\mathcal{M}:\gamma(0)=x, \gamma(1)=y \}.
\end{equation*} 
The curve(s) $\gamma$ achieving this minimum is called the shortest geodesic between $x$ and $y$.
However, the notion of geodesic distance between two points is not always obvious. In some cases, it may be useful to define the distance between two points on~$\mathcal{M}$ differently.

The gradient of a smooth scalar function $F:\mathcal{M}\rightarrow\Rbb$ at $x\in\mathcal{M}$ is the unique element $\grad_x F(x)\in \tgspace{x}{\mathcal{M}}$ that satisfies
\begin{equation*}
	\dirder{F}{x}{\xi} = g_x(\grad_x F(x),\xi), \quad \text{for all $\xi \in \tgspace{x}{\mathcal{M}}$},
\end{equation*}
where
\begin{equation*}
	\dirder{F}{x}{\eta} = \lim_{t \rightarrow 0} \frac{F(x+t\eta)-F(x)}{t}
\end{equation*}
is the standard directional derivative of $F$ at $x$ in the direction~$\eta$.

For quotient manifolds $\mathcal{M}=\overline{\mathcal{M}}/\sim$, where $\overline{\mathcal{M}}$ is the total space and $\sim$ is the equivalence relation that defines the quotient, the tangent space $\tgspace{\bar{x}}{\overline{\mathcal{M}}}$ at $\bar{x}$ admits a decomposition into its vertical and horizontal subspaces
\begin{equation*}
	\tgspace{\bar{x}}{\overline{\mathcal{M}}} = \mathcal{V}_{\bar{x}}  \oplus \mathcal{H}_{\bar{x}}.
\end{equation*}
The vertical space $\mathcal{V}_{\bar{x}}$ is the set of directions that contains tangent vectors to the equivalence classes. The horizontal space $\mathcal{H}_{\bar{x}}$ is a complement of $\mathcal{V}_{\bar{x}}$ in $\tgspace{\bar{x}}{\overline{\mathcal{M}}}$.
A tangent vector $\xi_x$ at $x\in\mathcal{M}$ has a unique representation $\bar{\xi}_{\bar{x}} \in \mathcal{H}_{\bar{x}}$ at $\bar{x}$. 
Provided that the metric $\bar{g}_{\bar{x}}$ in the total space is invariant along the equivalence classes, it defines a metric on the quotient space
\begin{equation*}
	g_{x}(\xi_{x},\zeta_{x}) \eqdef \bar{g}_{\bar{x}}(\bar{\xi}_{\bar{x}},\bar{\zeta}_{\bar{x}}).
\end{equation*}
If $\bar{F}$ is a function on $\overline{\mathcal{M}}$ that induces a function $F$ on $\mathcal{M}$, then
\begin{equation*}
	\overline{\grad_x F(x)} = \grad_{\bar{x}}\bar{F}(\bar{x}),
\end{equation*}
in which $\grad_{\bar{x}}\bar{F}(\bar{x})$ belongs to the horizontal subspace~$\mathcal{H}_{\bar{x}}$.

\clearpage

\addtocounter{sidebar}{1}
\section*{Sidebar \thesidebar: Basics Concepts of Local Sensitivity Analysis}

This briefly exposition recalls basic concepts of local sensitivity analysis. It follows the terminology of \cite{Khalil:2002wj}.

Consider an oscillator described by \eqref{eq:sys-para}.
Most characteristics of this system (defined in the previous sections) depend on the value of this parameter~$\lambda$. It means that, for each characteristic  of the system, there exists a function $c:\Lambda\rightarrow\mathcal{C}$ that associates with each value of the parameter $\lambda$ an element $c(\lambda)$ in the space $\mathcal{C}$ to which belongs the characteristic.

Under appropriate regularity assumptions (see \cite{Khalil:2002wj} for details), the sensitivity function~$S^c:\Lambda\rightarrow\tgspace{c(\lambda)}{\mathcal{C}}$ of the characteristic $c(\lambda)$ associates with each value of the parameter $\lambda$ the element~$S^c(\lambda)$ in the tangent space $\tgspace{c(\lambda)}{\mathcal{C}}$ at $c(\lambda)$, defined as
\begin{equation*}
	S^c(\lambda) \eqdef \frac{\partial c}{\partial \lambda}(\lambda) = \lim_{h\rightarrow 0}  \frac{c(\lambda+h) - c(\lambda)}{h}.
\end{equation*}
The sensitivity $S^c(\lambda)$ provides a first-order estimate of the effect of parameter variations on the characteristic. It can also be used to approximate the characteristic when $\lambda$ is sufficiently close to its nominal value $\nlambda$.
For small $\norm{\lambda - \nlambda}_2$, the characteristic $c(\lambda)$ can be expanded in a Taylor series about the nominal solution $c(\nlambda)$ to obtain
\begin{equation*}
	c(\lambda) = c(\nlambda) + S^{c}(\nlambda) \, \norm{\lambda - \nlambda}_2 + \mathcal{O}\left(\norm{\lambda - \nlambda}_2^2\right).
\end{equation*}
This means that the knowledge of the nominal characteristic $c(\nlambda)$ and the sensitivity function suffices to approximate the characteristic for all values of~$\lambda$ in a small ball centered at~$\nlambda$.

The main difficulty of sensitivity analysis is to formulate the appropriate (analytical) equation to be solved in order to find the characteristic $c(\lambda)$. Then, differentiating this (analytical) problem yields the sensitivity equation to be solved in order to find the sensitivity function $S^c(\nlambda)$. The analytical problem can be an algebraic problem, an initial value problem, a boundary value problem, etc.

\begin{rem}	
	If, for a given value of the parameter $\lambda$, the characteristic $c(\lambda)$ is itself a function $c(\lambda) : A \rightarrow B$ in the space of functions $\mathcal{C}$, the sensitivity $S^c(\lambda)$ is also a function $S^c(\lambda) : \tilde{A} \rightarrow \tilde{B}$ in the tangent space $\tgspace{c(\lambda)}{\mathcal{C}}$, where $\tilde{A}$ and $\tilde{B}$ are the domain and the image of the sensitivity function.
	For convenience, the characteristic and the sensitivity function are denoted by $c:A \times \Lambda \rightarrow B$ and $S^c : \tilde{A} \times \Lambda \rightarrow \tilde{B}$, respectively. 
\end{rem}

\begin{rem}
	It is often meaningful to compute the relative sensitivity function $\sigma^c(\lambda)$, defined as
	\begin{equation*}
		\sigma^c(\lambda) \eqdef \frac{\lambda}{\norm{c(\lambda)}_{c(\lambda)}} \, \frac{\partial c}{\partial \lambda} (\lambda)  = \lim_{h\rightarrow 0}  \frac{[c(\lambda+h) - c(\lambda)]/\norm{c(\lambda)}_{c(\lambda)}}{[\lambda + h - \lambda]/\lambda},
	\end{equation*}
	where $\norm{\cdot}_{c(\lambda)}$  denotes the norm induced by the Riemannian metric $g_{c(\lambda)}\left(\cdot,\cdot\right)$ at ${c(\lambda)}$. A relative sensitivity function measures the relative change in the model characteristic to a relative change in the parameter value.
\end{rem}

\clearpage

\addtocounter{sidebar}{1}
\section*{Sidebar \thesidebar: Numerical Tools} \label{sec:num_tools}

Several numerical algorithms exist for the numerical computation of periodic orbits \cite{Ascher:1995ty,Lust:2001iy,Seydel:2010tw}. Most algorithms recast the periodic orbit computation as a two-point boundary value problem. Numerical boundary value methods fall into two classes: 
\begin{enumerate}
	\item shooting methods generate trajectory segments using a numerical time integration and match segment endpoints with each other and the boundary conditions;
	\item global methods project the differential equations onto a finite dimensional space of discrete closed curves that satisfy the boundary conditions.
\end{enumerate}
Both methods yield a set of (nonlinear) equations that are solved with root-finding algorithms, usually Newton's method.

This sidebar summarizes popular algorithms for the computation of periodic orbits. Then it emphasizes how the computation of the infinitesimal phase response curve is a cheap by-product of this computation. Finally, it extends these algorithms for the computation of oscillator sensitivities: angular frequency, steady-state periodic solution, and infinitesimal phase response curve sensitivities. 
More sophisticated algorithms can be found in the literature and adapted similarly (see \cite{Guckenheimer:2000ea,Lust:2001iy,Govaerts:2006bv}).

\subsection{Numerical Computation of Periodic Orbits}

A periodic orbit $\gamma$ is characterized by the $2\pi$-periodic steady-state solution $\persol:\Sbb^1\rightarrow\gamma$ describing a closed curve in the state space and the angular frequency $\omega>0$ (or equivalently the period~$T$) that solve the boundary value problem \eqref{eq:bvp_x}.

Considering a (nonuniform) partition $\mesh$ of the unit circle $\Sbb^1$
\begin{equation} \label{eq:partition}
	\mesh \eqdef \{ 0 = \theta_0 < \theta_1 < \cdots < \theta_N = 2\pi \},
\end{equation} 
the $2\pi$-periodic steady-sate solution $\persol(\cdot)$ is numerically approximated by a closed discrete curve in the state space $\mathcal{X}$. A discrete curve is a set of points $\{\persol_0,\persol_1,\ldots,\persol_N\}$ associated with the set of phases \eqref{eq:partition}, such that $\persol_i$ approximates $\persol(\theta_i)$ for all $i=0,1,\ldots,N$. This discrete curve is closed, that is, $\persol_N = \persol_0$, which reflects the periodicity of the solution $\persol(\cdot)$. 
Below, the circle partition $\mesh$ is fixed and the discrete curve is numerically represented by the vector $\persol_\mesh \eqdef (\persol_0,\persol_1,\ldots,\persol_N)$.
Phase-steps are denoted by $h_i=\theta_{i+1} - \theta_i$.

Equations for approximate periodic orbits take then the form of $N$ $n$-dimensional vector equations
\begin{equation*}
	r_i(\persol_\mesh,\omega) = 0, \quad i=0,1,\ldots,N-1,
\end{equation*}
where different residual maps $r_i$ lead to different numerical methods (see \mytablename~\ref{tab:num_tools} for two popular one-step schemes).
These equations are completed by the periodicity condition
\begin{equation*}
	r_{N}(\persol_\mesh,\omega) \eqdef \persol_N - \persol_0 = 0
\end{equation*}
and the phase condition
\begin{equation*}
	r_{\PLC}(\persol_\mesh,\omega) \eqdef \PLC(\persol_\mesh;\lambda) = 0.
\end{equation*}

This set of (nonlinear) equations $r(\persol_\mesh,\omega) = 0$ is solved with the root-finding Newton's method. Starting from an initial guess $\left({(\persol_\mesh)}^{(0)},\omega^{(0)}\right)$, this method iteratively updates the solution
\begin{align*}
	{(\persol_\mesh)}^{(k+1)} & = {(\persol_\mesh)}^{(k)} + {(\Delta \persol_\mesh)}^{(k)}  \\
	\intertext{and}
	\omega^{(k+1)} & = \omega^{(k)} + \Delta \omega^{(k)}.
\end{align*}
Update terms are computed by solving the linear problem
\begin{equation} \label{eq:num_perorb}
	\begin{bmatrix}
		A       & b^{\persol} \\
		\trans{{c^{\persol}}} & d^{\persol} 
	\end{bmatrix}
	\begin{bmatrix}
		\Delta \persol_{\mesh} \\
		\Delta \omega
	\end{bmatrix}
	=
	-
	\begin{bmatrix}
		r_{\mesh}(\persol_\mesh,\omega) \\ r_{\PLC}(\persol_\mesh,\omega)
	\end{bmatrix},
\end{equation}
where $A$ has a particular block structure for one-step schemes and $b^{\persol}$, $c^{\persol}$, and $d^{\persol}$ are also defined by blocks
\begin{align*}
	A & = 
	\begin{bmatrix}
		 G_{0}	& -H_{0}	&        	&     		\\
				& \ddots 	& \ddots 	&     		\\
		      	&        	& G_{N-1}	& -H_{N-1} 	\\
		-I_{n}	&        	&        	&  I_{n}  
	\end{bmatrix},
	&
	b^{\persol} & =
	\begin{bmatrix}
		b^{\persol}_0    \\
		\vdots  \\
		b^{\persol}_{N-1} \\
		0_{n\times1}
	\end{bmatrix},
	\\
	\trans{{c^{\persol}}} & = 
	\begin{bmatrix}
		\frac{\partial \PLC}{\partial x_0} & \cdots & \frac{\partial \PLC}{\partial x_{N-1}} & \frac{\partial \PLC}{\partial x_{N}}
	\end{bmatrix},
	&
	d^{\persol} & = 
	\begin{bmatrix}
		\frac{\partial \PLC}{\partial \omega}
	\end{bmatrix}.
\end{align*}
Expressions for block entries $G_i$, $H_i$, and $b_i^{\persol}$ depend on the methods used to generate residual maps $r_i(\persol_\mesh,\omega)=0$, with $i=0,1,\ldots,N-1$, for approximate periodic orbits (see \mytablename~\ref{tab:num_tools}).

The main computational effort in one iteration is the evaluation of the  $(N+1)n \times (N+1)n$ structured matrix~$A$, whose block entries are computed through fundamental solution time integrations or Jacobian matrix evaluations.

\subsection{Numerical Computation of Phase Response Curves}

The infinitesimal phase response curve $\iPRCu:\Sbb^1\rightarrow\Rbb^n$ of a periodic orbit is calculated by applying \eqref{eq:dirder} that involves computing the gradient of the asymptotic phase map evaluated along the periodic orbit, that is, the function $p(\cdot)$.

The gradient of the asymptotic phase map evaluated along the periodic orbit $\iPRCx:\Sbb^1\rightarrow\Rbb^n$ is the solution of the boundary value problem \eqref{eq:bvp_q}.

The gradient is numerically approximated by a closed discrete curve, that is, a set of points $\{\iPRCx_0,\iPRCx_1,\ldots,\iPRCx_N\}$ associated with the set of phases~\eqref{eq:partition}, such that $\iPRCx_N = \iPRCx_0$. 
This discrete curve is numerically represented by the vector $\iPRCx_\mesh \eqdef (\iPRCx_0,\iPRCx_1,\ldots,\iPRCx_N)$.

Following the same procedure as for approximate periodic orbits, equations for approximate gradients take the form of $(N+1)n$ linear equations 
\begin{equation*}
	 \tilde{A} \, \iPRCx_\mesh = 0,
\end{equation*}
where the matrix $\tilde{A}$ has the same structure as the matrix $A$
\begin{equation*}
	\tilde{A} = 
	\begin{bmatrix}
		 \tilde{G}_{0}	& -\tilde{H}_{0}	&        			&     				\\
						& \ddots 			& \ddots 			&     				\\
		      			&        			& \tilde{G}_{N-1}	& -\tilde{H}_{N-1} 	\\
		-I_{n}			&        			&        			&  I_{n}  
	\end{bmatrix}.
\end{equation*}
Block entries of $\tilde{A}$ can be constructed based on numerical computations  for the periodic orbit computation (see \mytablename~\ref{tab:num_tools}).

The matrix $\tilde{A}$ is, by construction, singular with a simple rank deficiency. This rank deficiency is overcome by adding a normalization condition for $\iPRCx_\mesh$. Discretizing \eqref{eq:bvp_q_c} yields
\begin{equation*}
	\trans{v_\mesh} \, P \, \iPRCx_\mesh = \omega,
\end{equation*} 
where $v_\mesh\eqdef(f(\persol_0,0),f(\persol_1,0),\ldots,f(\persol_N,0))$ is the approximate tangent vector to the periodic orbit and $P$ is a ponderation matrix that depends on the method class. A standard method to obtain a system of defining equations that is square and regular is to border the matrix $\tilde{A}$ (see \cite[Theorem~5.8]{Seydel:2010tw} for details)
\begin{equation} \label{eq:num_iprc}
	\begin{bmatrix}
		\tilde{A} & b^\iPRCx \\
		\trans{{c^{\iPRCx}}} & d^\iPRCx 
	\end{bmatrix}
	\begin{bmatrix}
		\iPRCx_{\mesh} \\
		\xi
	\end{bmatrix}
	=
	\begin{bmatrix}
		0 \\ \omega
	\end{bmatrix},
\end{equation}
with $d^\iPRCx\neq0$, $\trans{{c^{\iPRCx}}} = \trans{v_\mesh}P$, and $b^\iPRCx\notin\range(\tilde{A})$ (for example $b^\iPRCx=v_\mesh$).

\subsection{Numerical Computation of Oscillator Sensitivities}

The angular frequency sensitivity $S^\omega\in\Rbb^{1\times l}$ and the  sensitivity of the $2\pi$-periodic steady-sate solution  $S^{\persol}:\Sbb^1\rightarrow\Rbb^{n\times l}$ are the solutions of the linear boundary value problem \eqref{eq:bvp_Zx}. 
Equations for approximate periodic orbit sensitivities take the form of a system of linear equations
\begin{equation} \label{eq:num_sens_perorb}
	\begin{bmatrix}
		A      & b^{\persol} \\
		\trans{{c^{\persol}}} & d^{\persol}
	\end{bmatrix}
	\begin{bmatrix}
		S^{\persol}_\mesh \\
		S^\omega
	\end{bmatrix}
	=
	\begin{bmatrix}
		E^{\persol}_{\mesh} \\ E^{\persol}_{\PLC}
	\end{bmatrix},
\end{equation}
where $E^{\persol}_{\PLC} = -\frac{\partial \PLC}{\partial \lambda}$ and $E^{\persol}_i$ depends on the numerical method used (see \mytablename~\ref{tab:num_tools}).

The sensitivity of the gradient of the asymptotic phase map evaluated along the periodic orbit $S^\iPRCx:\Sbb^1\rightarrow\Rbb^{n\times l}$ is the solution of the linear boundary value problem \eqref{eq:bvp_Z_iPRCx}. 
Equations for approximate infinitesimal phase response curve sensitivities take the form of a system of linear equations
\begin{equation} \label{eq:num_sens_iprc}
	\begin{bmatrix}
		\tilde{A} & b^{\iPRCx} \\
		\trans{{c^{\iPRCx}}} & d^{\iPRCx} 
	\end{bmatrix}
	\begin{bmatrix}
		S^{\iPRCx}_{\mesh} \\
		\xi
	\end{bmatrix}
	=
	\begin{bmatrix}
		E^{\iPRCx}_{\mesh} \\ E^{\iPRCx}_\omega
	\end{bmatrix},
\end{equation}
where $E^\iPRCx_\omega = S^\omega - \trans{{S^{v}_{\mesh}}}\,P\,\iPRCx_\mesh$ and $E^{\iPRCx}_i$ depends on the numerical method used (see \mytablename~\ref{tab:num_tools}).

In \eqref{eq:num_sens_perorb} and \eqref{eq:num_sens_iprc}, the square matrices in left-hand sides are identical to the matrices used for the computation of the periodic orbit in \eqref{eq:num_perorb} and the gradient in \eqref{eq:num_iprc}, respectively. The only additional computation effort arises from the evaluation of the right-hand sides.

\clearpage

\section*{Author Information}

\noindent Pierre Sacr\'e (S'10) received the M.S.~degree in aerospace engineering (2008) and the Ph.D.~degree in engineering sciences (2013), both from the Universit\'e de Li\`ege, Belgium.
	
He was a Research Fellow of the Belgian National Research Fund (F.R.S.-FNRS) in the Department of Electrical Engineering and Computer Science at the Universit\'e de Li\`ege, Belgium.
He was a visiting postdoctoral research associate in the Department of Mathematics at the Imperial College London in Fall 2013. Since January 2014, he has been a Fulbright Belgium research scholar in the Department of Biomedical Engineering at Johns Hopkins University, Baltimore, Maryland. His research interests include dynamical systems, oscillators, and biological applications.

\bigskip \bigskip \bigskip \bigskip

\noindent Rodolphe Sepulchre (F'10) received the engineering degree (1990) and the Ph.D. degree (1994), both in mathematical engineering, from the Universit\'e catholique de Louvain, Belgium.

He was a BAEF Fellow in 1994 and held a Postdoctoral Position at the University of California, Santa Barbara from 1994 to 1996. He was a Research Associate of the FNRS at the Universit\'e catholique de Louvain from 1995 to 1997. Since 1997, he has been Professor in the Department of Electrical Engineering and Computer Science at the Universit\'e de Li\`ege. He was Department Chair from 2009 to 2011. He held Visiting Positions at Princeton University (2002--2003) and Mines Paris-Tech (2009--2010) and Part-Time Positions at the Universit\'e catholique de Louvain (2000--2011) and at INRIA Lille Europe (2012--2013). Since 2013, he is Professor of Engineering at the University of Cambridge.

In 2008, he was awarded the IEEE Control Systems Society Antonio Ruberti Young Researcher Prize. He is an IEEE Fellow and an IEEE CSS Distinguished Lecturer since 2010.

\bigskip \bigskip 

Email: r.sepulchre@eng.cam.ac.uk

Mailing address: University of Cambridge, Department of Engineering, Trumpington Street, Cambridge CB2 1PZ, United Kingdom



\nomenclature[aa]{$u$}{input value of a system}
\nomenclature[ab]{$y$}{output value of a system}
\nomenclature[ac]{$x$}{state variable of a system}

\nomenclature[ad]{$\perorb$}{periodic orbit of an oscillator}
\nomenclature[ae]{$\persol$}{zero-input steady-state solution of an oscillator}

\nomenclature[ae]{$\Theta$}{asymptotic phase map of an oscillator}
\nomenclature[af]{$\theta$}{phase variable of an oscillator}

\nomenclature[bb]{$Q$}{(finite) phase response curve}
\nomenclature[bc]{$q$}{infinitesimal phase response curve}

\nomenclature[bd]{$p$}{gradient of the asymptotic phase map evaluated along the periodic orbit, $\iPRCx(\cdot) \eqdef \nabla_{x}\Theta(\persol(\cdot))$}

\nomenclature[ca]{$\Rbb$}{set of real numbers}
\nomenclature[cb]{$\Rbb^n$}{$n$-dimensional Euclidean space}
\nomenclature[cc]{$\Sbb^1$}{set of points on the unit circle, $\Sbb^1 \eqdef \Rbb/(2\pi\Zbb)$}
\nomenclature[cd]{$\Nbb$}{set of natural numbers}

\nomenclature[da]{$\dot{x}$}{derivative, with respect to time, of the variable $x$}
\nomenclature[db]{$x'$}{derivative, with respect to phase, of the variable $x$}
\nomenclature[dc]{$\zeroinput$}{input identically equal to $0$}

\nomenclature[ma]{$\cconj{z}$}{complex conjugate of the complex number $z$}
\nomenclature[mb]{$\trans{A}$}{transpose of the matrix $A$}

\nomenclature[q]{$(x,y)$}{equivalent notation for the vector $\trans{[\trans{x} \trans{y}]}$}

\nomenclature[q]{$\norm{x}_{\infty}$}{maximum norm of the vector $x$, $\norm{x}_{\infty}\eqdef \max(\abs{x_1},\ldots,\abs{x_n})$}

\end{document}